\documentclass[a4paper,10pt]{amsart}%{article}

\usepackage{color}
\usepackage{amssymb}
\usepackage[pdftex, colorlinks=false, backref]{hyperref}

\usepackage{graphicx}
\usepackage{psfrag}

\newcommand{\n}{\noindent}
\newcommand{\R}{\mathbb{R}}
\newcommand{\N}{\mathbb{N}}

\newcommand{\Lip}{\mathrm{Lip}}

\newcommand{\BV}{\mathbf{BV}}
\newcommand{\ve}{\varepsilon}
\newcommand{\Lsp}{\mathbf{L}}
\newcommand{\con}{\mathbf{C}}
\newcommand{\dott}{\,\cdot\,}
\newcommand{\esssup}{\mathrm{ess}\sup}

\newtheorem{theorem}{Theorem}[section]

\newtheorem{lemma}[theorem]{Lemma}
\newtheorem{definition}[theorem]{Definition}
\newtheorem{example}[theorem]{Example}
\newtheorem{remark}[theorem]{Remark}

\numberwithin{equation}{section}     
\allowdisplaybreaks

\title{On an inverse problem for scalar conservation laws}

\begin{document}

\author[Holden]{Helge Holden} \address[Holden]{\newline Department of
  Mathematical Sciences, Norwegian University of Science and
  Technology, NO--7491 Trondheim, Norway,\newline {\rm and} \newline
  Centre of Mathematics for Applications,
  % Department of Mathematics,
  University of Oslo, P.O.\ Box 1053, Blindern, NO--0316 Oslo, Norway
} \email[]{holden@math.ntnu.no} \urladdr{www.math.ntnu.no/\~{}holden}

\author[Priuli]{Fabio Simone Priuli} \address[Priuli] {\newline
  Dipartimento di Matematica, Universit\`a degli Studi di Roma Tor Vergata, 
  Via della Ricerca Scientifica 1, I--00133 Roma, Italy.}
\email[]{priuli@mat.uniroma2.it}
\urladdr{http://www.mat.uniroma2.it/\~{}priuli/}

\author[Risebro]{Nils Henrik Risebro} \address[Risebro]{\newline
  Centre of Mathematics for Applications,
  % Department of Mathematics,
  University of Oslo, P.O.\ Box 1053, Blindern, NO--0316 Oslo,
  Norway.}  \email[]{nilshr@math.uio.no}
\urladdr{www.math.uio.no/\~{}nilshr}

\date{\today} 

\thanks{Supported in part by the Research Council of Norway.}

\begin{abstract}
  We study in what sense one can determine the flux functions $k=k(x)$ and $f=f(u)$,
  $k$ piecewise constant, in the scalar hyperbolic conservation law 
  $u_t+(k(x)f(u))_x=0$ by observing the solution $u(t,\dott)$ of the Cauchy problem with
  suitable piecewise constant initial data $u|_{t=0}=u_o$.
\end{abstract}

\maketitle

\section{Introduction}

In this paper, we deal with the inverse problem for scalar
conservation laws. More precisely, we consider a scalar conservation
law of the form
\begin{equation}\label{eq:intrononh}
  \partial_t u + \partial_x \big(k(x)f(u)\big) =0\,,
\end{equation}
with $(t,x)\in[0,\infty)\times\R$, $u(t,x)\in\R$,
$k\colon\R\to\,(0,\infty)$ and $f\colon\Omega\subseteq\R\to\R$ flux
functions whose smoothness will be prescribed later. It is well known
that if $k$ is a constant function and $f$ is locally Lipschitz
continuous, then for every initial data
$u_o\in\Lsp^\infty(\R)\cap\Lsp^1(\R)$ there exists a unique entropy
solution $u(t,\dott)\in\con^0([0,\infty),\Lsp^1(\R))$,
see~\cite{BressanBook, DafermosBook, HoldenRisebroBook}. In recent
years, motivated by problems arising in traffic flow
models~\cite{KlausenRisebroJDE,KlingRisebro, Towers} and in multiphase
flow models in porous media~\cite{Bachmann,Diehl, Diehl2,GimseRisebro}, the equation \eqref{eq:intrononh} has been widely studied also in the
case where $k$ is a discontinuous and piecewise constant function. In
this latter case, assuming that the flux function $f$ is
strictly concave and defined in a compact interval $[u_1,u_2]$ with
$f(u_1)=f(u_2)=0$, it has been proved in~\cite{KlingRisebro} that a
unique entropy solution exists for every initial data in $\BV(\R)$.

\medskip

The goal of this paper is to find a reconstruction procedure which
allows us to approximate the unknown functions $k$ and $f$
in~\eqref{eq:intrononh} starting from the observation of the solutions
$u_{\rm obs}(t,x)$ corresponding to Cauchy problems with suitably
chosen initial data. This is a so-called \emph{coefficient inverse
  problem}, because an observer has complete access to both initial
data and solutions of the problem, but only partial information on
the structure of the equation itself. The results of our work ensure
that the reconstruction is possible for two important classes of problems:
when $k\equiv const$ and $f$ is  sufficiently smooth, and when $k$ is
piecewise constant and $f$ is a known strictly concave function.

\medskip

This kind of inverse problem has many applications, depending on the
underlying physical phenomena described by~\eqref{eq:intrononh}. For
instance, we can consider models of traffic flow on highways
(see~\cite{KlingRisebro,LW,Richards}). Here the unknown $u(t,x)$
denotes the density of cars at time $t$ in the position $x$, the
product $k(x)f(u(t,x))$ represents the flux of cars which cross each
position $x$ at a time $t$
per unit of time, and the function $k(x)$ describes specific
characteristics of the road in the position $x$. The inverse problem,
in this case, corresponds to the problem of determining the unknown
properties $k$ and $f$ of the considered road by only monitoring the
resulting density of cars $u_{\rm obs}$ along the road.

Also, we want a procedure that can handle problems where parts of the
spatial domain are not directly observable and, hence, where data from
the observable regions has to be used to reconstruct the
characteristics of the physical system also in the unobservable
regions.  To fix the ideas, think of a highway where a tunnel is
present in an interval $[a,b]$, or where the traffic data is monitored
by using sensors which cannot cover the whole road. In this situation,
in addition to reconstructing $k$ and $f$ in the observable region
$\R\setminus [a,b]$, we would like to determine the flux function
$k(x)$ in $[a,b]$, relying on the observed data $u_{\rm
  obs}|_{\R\setminus [a,b]}$, to detect the possible obstructions, due
to car accidents or other events, and to locate their precise position
inside the region $[a,b]$.

% \medskip

Despite the ample spectrum of applications, to our knowledge only few
attempts of addressing inverse problems for conservation
laws~\eqref{eq:intrononh} have been made.

In~\cite{KangTanuma} a special class of inverse problems is solved for
scalar conservation laws~\eqref{eq:intrononh} with $k(x)\equiv 1$ and
$f$ of class $\con^2$ and uniformly convex. Namely, it is assumed that
the initial data for~\eqref{eq:intrononh} is such that the observed
solution $u_{\rm obs}$ consists only of a single shock wave, after a
large enough time $T$. In this particular case, $f$ can be expressed as
limit of functions explicitly depending on the shock wave and on the
initial data. Unfortunately, the requirement that the solution
develops a single discontinuity is very strong in the context of
conservation laws, making this approach infeasible for general
equations of the form~\eqref{eq:intrononh}.

In~\cite{JamesSepulvedaSIAM} and then in~\cite{BBCS,CJS,FBS}, a more general approach is presented to
deal with the inverse problem for~\eqref{eq:intrononh} under the
assumptions of $k(x)\equiv 1$ and $f$ locally Lipschitz
continuous. Namely, the flux function $f$ is uniquely identified by
minimizing, over a compact set of Lipschitz continuous fluxes, a
suitable cost functional $J(f)$ which measures the distance between
the observed solution $u_{\rm obs}$ and the solution corresponding to
any choice of the flux. The functional has the following form
\begin{equation}\label{eq:JScost}
  J(f):= \,{1\over 2}\,\|u_f (T,\dott)-u_{\rm obs}\|_{\Lsp^2}^2 +\,{\rho\over 2}\, \left| \int_\R x(u_f(T,x)-u_{\rm obs}(x))\, dx \right|
\end{equation}
where $u_f$ is the solution to the conservation law with flux $f$ and
$\rho>0$ is a fixed constant. The first term is the cost used in the
well-known output least square method and is sensitive to the shape of
the observed function, while the second term is more sensitive to the
localization of the observed function on the $x$-axis. Dealing with a
minimization problem for~\eqref{eq:JScost}, differentiability of $J$
with respect to $f$ is important, since both optimality conditions and
gradient algorithms rely on it, however, in general the function is
nondifferentiable. Yet, minimization is possible if additional
assumptions are posed on the number and location of jumps in the
observed solution $u_{\rm obs}$. Unfortunately, one cannot in general
expect these additional hypotheses to hold, and in the general setting
the problem remains open.

Similar results are obtained in \cite{CastroZuazua} where the flux
$f$ is obtained by minimizing the functional
\begin{equation*}
  J(f):= \frac12\| u_f(T,\dott) - u_{\rm obs}\|_{\Lsp^2}^2 +
  \frac{\rho}{2}
  \int_{u_1}^{u_2} | f'(u)| \,du.
\end{equation*}
If the penalization parameter $\rho$ is zero, then the above
functional does not have a unique minimizer as can easily be
demonstrated by an example where $u_{\rm obs}$ contains
shocks. Nevertheless, in \cite{CastroZuazua} efficient algorithms are
developed for the numerical calculation of minimizers $f$ even if the
observed solution has discontinuities.

Finally, the recent paper~\cite{BD} develops an alternative reconstruction 
method, still based on a constrained minimization procedure, for a specific 
model~\eqref{eq:intrononh} with $k\equiv 1$ and $f$ of class $\con^2$ 
with a single inflection point, describing the sedimentation of small particles dispersed
in a viscous fluid.

In this paper we follow a different approach. We exploit the complete
and detailed knowledge of the approximation procedure used to obtain
solutions to the Cauchy problem for~\eqref{eq:intrononh}, the
so-called \emph{front-tracking} algorithm~\cite{BressanBook,
  HoldenRisebroBook}, in order to somehow revert the construction and
deduce properties of the flux functions $k,f$ starting from the
observed solutions. Our analysis is restricted to one space dimension due to the constructive method that we advocate. For applications to traffic flow, this suffices.

This produces an \emph{ad hoc} procedure which allows us to solve the
inverse problem, both in the case of homogeneous conservation laws
where $k\equiv \text{const}$, and in the case of a piecewise constant
function $k$, as long as we assume that we can observe the solutions
corresponding to suitable families of initial data.  Namely, for the
case $k\equiv 1$, i.e., for the homogeneous conservation law
\begin{equation}\label{eq:introhom}
  \partial_t u + \partial_x f(u) =0\,,
\end{equation}
we prove in Theorem~\ref{thm:1} the following: If $f$ is of class
$\con^{1,1}$ with a finite number of inflection points, then we can
always find a piecewise affine interpolation $f_\nu$ of $f$, by using
a single observation at a fixed time $T>0$ of a finite number of
solutions $u_{\rm obs}$, corresponding to properly chosen initial
data. Such approximate flux $f_\nu$ coincides with $f$ in suitable
nodes $u_1<\dots<u_\nu$, and it is close to $f$ in the sense that the
$\Lsp^1$ distance between $u_{\rm obs}$ and the solution of the
conservation law with flux $f_\nu$ converges to $0$ as $\nu\to\infty$.

To deal with the general case of a piecewise constant function $k(x)$,
we focus our attention on the case when $f(u)$ is a known function and only
$k$ has to be reconstructed. This assumption does not seem to be 
unnatural since we  can expect there exists some observable spatial 
interval $J$ where $k$ is known to be constant, say $k(x)\equiv k_o$, and 
in such a region Theorem~\ref{thm:1} allows to reconstruct a good approximation
of the function $f$. In the traffic flow example, you can think to first reconstruct
$f$ by observing the car behavior in a small portion of road that you know it
is spatially homogeneous, and then to use this knowledge to reconstruct the
inhomogeneities in the rest of the road. 

In Theorem~\ref{thm:2}, assuming that $f$ is defined on an interval
$[u_1,u_2]$, is strictly concave and such that $f(u_1)=f(u_2)=0$
(which is the case, e.g., in the Lighthill--Whitham--Richards traffic
flow model~\cite{LW, Richards}), we prove that in order to reconstruct
exactly the function $k(x)$ on any compact interval $J\subseteq \R$,
it is enough to observe the solution $u_{\rm obs}$ in $[0,T]\times\R$,
% i.e., for all times in $[0,T]$ and in the whole spatial domain,
for a single suitable initial data $u_o^J$.

Finally, under the same assumptions on $f$, we have studied the case 
in which the solution can only be
observed in $[0,T]\times (\R\setminus I)$, for some unobservable open
interval $I$ and for some time $T$ large enough. In this case, the
expression of $k(x)$ outside $I$ can be obtained by using
Theorem~\ref{thm:2}, but $k(x)$ can also be reconstructed inside $I$,
if we assume that no more than two jumps are present inside the
unobservable interval. Namely, in Theorem~\ref{thm:3} we prove that a
suitable choice of the initial data in the region $\{x\in\R~;~x<\inf
I\}$ allows us to reconstruct the position and the size of the jumps
of $k(x)$ inside $I$ from the observed solution $u_{\rm
  obs}$. Moreover, in Theorem~\ref{thm:4} we prove that the
reconstruction is also possible when the initial data cannot be chosen
freely but it is given by a constant state $\bar u_o$. This is for
instance the case when considering a physical system whose
inhomogeneity appears at time $t=0$, due to some external event (like
a car accident) which modifies the properties of the flux function in
a specific region. In this latter case, we prove that it is still
possible to determine positions and sizes of the jumps of $k(x)$ in
$I$, provided that the jump is large enough to influence the dynamics
outside $I$.

We remark that the assumption on the number of jumps in the
unobservable region $I$ is rather strong, because we are basically assuming
that only a single obstruction can be present. However, this appears to be unavoidable,
because if more than two jumps are allowed in $I$, then the inverse
problem is in general ill-posed. Indeed, in
Section~\ref{sec:counterex} we present a few examples where relaxing
the assumption on $k$ leads to infinitely many piecewise constant
functions $k(x)$ on $I$,  all giving the same observed solution in
$[0,T]\times (\R\setminus I)$. This means that in many real situations it is impossible, based only on the observations of the
solution in $[0,T]\times (\R\setminus I)$,
to distinguish between 
a single large obstruction or many smaller ones, . In such a context, one can 
apply Theorem~\ref{thm:3} in order to obtain a reconstructed flux function $k$ with a single jump
and consider such a single obstruction as an approximation of the real one, whose
structure can be very complex.

This paper represents the first steps towards a more complete understanding of the \emph{coefficient inverse problem}. Further study is necessary in order to address the fundamental question of stability.  Furthermore, extensions to multi-dimensional cases, will require novel techniques. 

\section{Main results}\label{sec:main_res}
We start by studying the inverse problem for~\eqref{eq:introhom},
i.e., in the case of $k(x)\equiv \text{const}$. We recall that a
Riemann problem for~\eqref{eq:introhom} is a Cauchy problem with
initial data of the form
\begin{equation}\label{eq:riem_data}
  u_o(x)=\left\{
    \begin{array}{ll}
      u^-, & \,\,x<0,\\
      u^+, & \,\,x>0,
    \end{array}
  \right.
\end{equation}
for given values $u^-\ne u^+$. 
In the following, it is fundamental to specify in which sense we observe the solution
to~\eqref{eq:intrononh} or~\eqref{eq:introhom}, and to this purpose we precisely
introduce next definition.  Please note that the use of ``observable'' 
in this paper differs from that in control theory. 

\begin{definition} \label{def:observable}
 {\it (a)} A function $z\colon \R\to\R$ is said to be \emph{observable} if we know its
   values $z(x)$ for (almost) every $x\in\R$.\\
 {\it (b)} Fixed $T>0$ and an interval $[a,b]\subset \R$, a function 
   $z\colon[0,\infty)\times\R\to\R$ is said to be \emph{partially observable} if
   $z(t,\cdot)\bigm|_{\R\setminus(a,b)}$ is observable for all $t\in [0,T]$ 
   in the sense of {\it (a)}.
\end{definition}

The choice to require observations of $z$ in the whole $\R$ (resp. 
$\R\setminus(a,b)$) has been made for sake of simplicity. The proof of 
most results presented in this paper could be adapted to the case of 
observability of $z$ on given bounded intervals, provided the model 
under consideration justify some a priori bound to the propagation 
speed (or equivalently to $f'$).

Our first result states that if the flux function $f$ is piecewise smooth and it has a finite number of inflection points
and if all solutions $u_{\rm obs}(T,\dott)$ to
Riemann problems~\eqref{eq:introhom}--\eqref{eq:riem_data} at some
fixed time $T>0$ are observable, then it is possible to
construct on every bounded interval $I\subseteq \R$ a piecewise
affine interpolation $\tilde f$ of the flux $f$, which is close to $f$
in the following sense: at a time $T$, the solution to every Cauchy
problem for
\begin{equation*}
  \partial_t u + \partial_x \tilde f(u) =0\,,
\end{equation*}
is close in $\Lsp^1$ to the solution to~\eqref{eq:introhom} with the
same initial data. More precisely, we prove the following:

\begin{theorem}\label{thm:1} 
  Let $T>0$, $u_*,u^*\in\R$ such that $u_*<u^*$
  and $c\in\R$ be fixed. Assume that $f\colon \Omega\to\R$ is continuous and piecewise
  $C^1$ with a finite number of inflection points on any bounded interval contained in $\Omega$,
  that $[u_*,u^*]\subseteq\Omega$, that $f(u_*)=c$ and that the solution to
  any Riemann problem for~\eqref{eq:introhom} at time
  $T$ is observable, in the sense of Definition~\ref{def:observable}{\it (a)}.  Then, for all $\nu\in\N$, setting
  $\delta:=2^{-\nu}|u^*-u_*|$ and $u_\alpha= u_*+\alpha\delta$ for
  $\alpha=0,\ldots,2^\nu$, there exists a piecewise affine function
  $f^\nu\colon [u_*,u^*]\to\R$ such that $f_\nu(u_\alpha)=f(u_\alpha)$
  for all $\alpha$ and
  \begin{equation}\label{eq:recest}
    \esssup_{[u_*,u^*]} |f'_\nu-f'| \le \Lip(f') \delta\,,
  \end{equation}
  where $\Lip(f')$ is a Lipschitz constant for the derivative $f'$ on
  $[u_*,u^*]$.

  This function $f_\nu$ represents a good reconstruction of the
  unknown flux $f$ in the following sense: if $\hat u$ is a $\BV$
  function with values in $[u_*,u^*]$, and we denote by $u^\nu$
  (resp.~$u_{\rm obs}$) the solution to the Cauchy problem for
  $\partial_t u + \partial_x f^\nu(u) =0$
  (resp.~for~\eqref{eq:introhom}) with initial data $\hat u$, then
  \begin{equation}\label{eq:close}
    \| u^\nu(T,\dott) - u_{\rm obs}(T,\dott) \|_{\Lsp^1(\R)} \le
    C T  \delta
  \end{equation}
  for a constant $C$ which does not depend on $\delta$.
\end{theorem}
We remark that~\eqref{eq:close} follows immediately
from~\eqref{eq:recest} and from the general stability results
contained in~\cite{HoldenRisebroBook} (see
Theorem~\ref{thm:lip_est}). Here the relevant result is the procedure
to construct a piecewise affine interpolation $f^\nu$ from the
observed solutions, so that $f^\nu$ coincides with the original flux
$f$ at points $u_o=u_*<u_1<\dots<u_{2^\nu}=u^*$ of the interval
$[u_*,u^*]$ and satisfies~\eqref{eq:recest}. Such a procedure is explicitly 
presented in Section~\ref{sec:proof} and constitutes the main part of the proof
of the theorem. We also stress that no
assumptions are made in Theorem~\ref{thm:1} concerning the regularity
of the observed solutions or concerning their discontinuity
structure. Furthermore, general solutions containing any finite number
of shocks and centered rarefaction waves can appear without affecting
the result of the reconstruction.

\begin{remark}\label{rem:inflections}
 Concerning the assumption on the finite number of inflection points of $f$ 
 in Theorem~\ref{thm:1}, it is important to notice that we are not prescribing 
 any knowledge of the actual location of the inflection points. This means 
 that for a fixed $\nu\in\N$ we have no idea of how close is $f_\nu$ to $f$: 
 if all the inflection points of $f$ are contained in the interior of a single interval 
 $I_\alpha =\,]u_\alpha,u_{\alpha+1}[$, $\alpha=0,\ldots,2^\nu$, 
 our reconstruction $f_\nu$ would be a monotone function with a completely 
 different behavior in $I_\alpha$, and only the estimate~\eqref{eq:recest} 
 would hold. What Theorem~\ref{thm:1} ensures is that for all $\nu\in\N$ the solution 
 $u^\nu$ satisfies~\eqref{eq:close} at time $T$ and that there exists $\nu$ large enough 
 so that at most one inflection point belongs to each interval $]u_\alpha,u_{\alpha+1}[$, 
 but we have no way to estimate a priori how large this $\nu$ must be.
 
 If modeling considerations could justify a lower bound on the distance between 
 consecutive inflection points, then it is immediate to verify that, choosing 
 $\delta$ smaller than this lower bound, we obtain a reconstruction $f_\nu$
 which captures much better the real shape of $f$.
 However, in the general case, there is no analogous strategy to apply and the
 best one can hope is to design system--specific adjustments in the 
 reconstruction, as we point out in Remark~\ref{rem:homog_improved}.
\end{remark}

Next we study scalar conservation laws of the more general
form~\eqref{eq:intrononh} with a piecewise constant term $k(x)$. Since
in general the existence of solutions to the Cauchy problem
for~\eqref{eq:intrononh} is much more difficult to prove than
for~\eqref{eq:introhom} (see, e.g., \cite{KarlsenRisebroTowers} and
references therein), we focus our attention on a specific class of
conservation laws studied in~\cite{KlausenRisebroJDE, KlingRisebro,
  Risebrointro,Towers} for which existence of a solution to the
Cauchy problem has been proved by Klingenberg and
Risebro~\cite{KlingRisebro}. Namely, we assume:
\begin{itemize}
\item[{\bf (H1)}] $\bullet$~$k\colon \R\to (0,\infty)$ is piecewise
  constant and belongs to $\BV(\R)$;
\item[{\phantom{\bf (H1)}}] $\bullet$~$f\colon[u_1,u_2]\to [0,\infty)$
  is of class $\con^2$, strictly concave and such that
  $f(u_1)=f(u_2)=0$. In particular, $f>0$ in $(u_1,u_2)$ and there
  exists a unique $u^m\in\,(u_1,u_2)$ such that
  $f(u^m)=\max_{[u_1,u_2]} f$.
\end{itemize}

\n From~\cite{KlingRisebro,Risebrointro}, we know that every Cauchy
problem for~\eqref{eq:intrononh} with flux functions $k$, $f$
satisfying {\bf (H1)}, and initial data in $\BV(\R)$, admits a unique
entropy solution in $\con([0,T];\Lsp^1(\R))$, see
Theorem~\ref{thm:KR}.

\begin{example}[{\bf Traffic flow on highways}]\label{ex:highway} A typical
  example of a system satisfying {\bf (H1)} is the simple inhomogeneous
  variant of the classical Lighthill--Whitham--Richards
  model~\cite{LW,Richards} for car traffic flow on a highway, obtained
  by multiplying the flux function $f(u)=u(1-u)$ with a piecewise
  constant factor $k(x)$. In this model, $u$ represents the density of
  cars on the highway and takes values in $[u_1,u_2]=[0,1]$ and $f(u)$
  represents the flux of cars per unit of time. The function $k(x)$
  represents specific features of the road considered in different
  spatial regions, e.g., regions in which cars have to reduce their
  speed or are allowed to increase it, all due to external factors.
\end{example}
Motivated by Example~\ref{ex:highway} above, in the following we will
say that a spatial region $I\subseteq \R$ is \emph{congested}
(resp.~\emph{fully congested}) if $u(x)\ge u^m$ (resp.~$u(x)\equiv
u_2$) for all $x\in I$.

We notice that if the flux functions satisfy assumptions {\bf (H1)}
and the solutions to any Riemann problem for~\eqref{eq:intrononh} are
observable, then we can first consider a small region $[\alpha,\beta]$ 
where the road is homogeneous and use Theorem~\ref{thm:1} on the 
interval $[u_1,u_2]$, with Riemann data centered in $x=(\alpha+\beta)/2$, 
to reconstruct $f(u)$ with a given precision. Hence, without loss of 
generality, we assume $f(u)$ to be a given function and focus our 
attention on the piecewise constant function $k(x)$. Under these 
assumptions, we can prove the following result, which provides an 
exact reconstruction procedure for the function $k(x)$ on any compact 
interval (see the proof of the theorem in Section~\ref{sec:proof}).

\begin{theorem}\label{thm:2} Let $T>0$ and $J\subseteq \R$ be a fixed
  compact interval. Assume that the function $f$
  in~\eqref{eq:intrononh} satisfies {\bf (H1)}, and that the solution to any
  Riemann problem for~\eqref{eq:intrononh} is observable for all times
  $t\in\,(0,T]$, in the sense of Definition~\ref{def:observable}{\it (a)}.  Then, there exists a
  unique piecewise constant function $k^J\colon J\to\R$ such that the
  following property holds. If we denote by $u^J$ (resp.~$u_{\rm
    obs}$) the solution to the Cauchy problem for $\partial_t u
  + \partial_x\big( k^J(x)f(u)\big) =0$
  (resp.~for~\eqref{eq:intrononh}) with initial data $\hat
  u\in\BV(\R)$ taking values in $[u_1,u_2]$, then
  \begin{equation}\label{eq:reconstr_2}
    u^J(t,x)=u_{\rm obs}(t,x)\,,\quad x\in J,\, t\in [0,T]\,.
  \end{equation}
\end{theorem}
Here, solutions must be observed on some interval $t\in\,(0,T]$ and
not only at a single time $t=T$. The reason for this additional
requirement is that, since the locations of the jumps in $k$ are
unknown, it is otherwise difficult to observe the speed of the waves
appearing in the solution. However, the time interval $(0,T]$ can be
taken arbitrarily small without interfering with our reconstruction
procedure.

Finally, we focus our attention to the case of incomplete
observability, i.e., when a part of the domain cannot be directly
observed. To fix the ideas, we assume that such unobservable part is
a given interval $(a,b)$.
Since Theorem~\ref{thm:2} can be used to reconstruct $k(x)$ on every
compact interval $J\subseteq\,(-\infty,a]$ and $J\subseteq[b,\infty)$,
it is not restrictive to assume that $k(x)$ is known and constant in
the observable region $\R\setminus (a,b)$. Moreover, we assume that in
the unobservable region $[a,b]$ the changes in $k$ can only be due to
some sort of obstruction which reduces the speeds of propagation. In
other words, we assume in the following that
\begin{itemize}
\item[{\bf (H2)}] $k(x)\equiv k_o$ in $\R\setminus(a,b)$ and $k(x)\le
  k_o$ for all $x\in[a,b]$.
\end{itemize}
To prove our main results for the inverse problem with partial
observability, we need to introduce a further hypothesis on the
function $k(x)$ in the unobservable interval $[a,b]$.
\begin{itemize}
\item[{\bf (H3)}] $k(x)$ has exactly two jumps in $[a,b]$, i.e., there
  exist $k_1\in\,(0,k_o)$ and $a\le\xi_1<\xi_2\le b$ such that
  \begin{equation}\label{eq:disck}
    k(x):= \begin{cases}
      k_o & x\notin (\xi_1,\xi_2)\,, \\
      k_1& x\in\,(\xi_1,\xi_2)\,.
    \end{cases}
  \end{equation}
\end{itemize}
By applying Theorem~\ref{thm:KR} we know that to each choice
$(k_1,\xi_1,\xi_2)$ in $(0,k_o)\times[a,b]\times[a,b]$ there
corresponds a flux function $k(x)$, defined by~\eqref{eq:disck}, such
that any Cauchy problem for~\eqref{eq:intrononh} with $\BV$ initial
data has a unique entropy solution.

For a scalar conservation law~\eqref{eq:intrononh} satisfying
hypotheses {\bf (H1)}--{\bf (H3)}, we consider two different inverse
problems, corresponding to two possible applications to the traffic
flow model described in Example~\ref{ex:highway}: Reconstruction from
initial data which is a stationary solution in $[a,\infty)$ and
reconstruction from a constant initial data.

The first problem is the reconstruction of $k(x)$ in the case of an
initial data $u(0,x)|_{[a,\infty)}=u_o(x)$ which is a stationary
entropy solution of~\eqref{eq:intrononh} in $[a,\infty)$ with values
in $[u_1,u_2]$. In other words, we assume that the initial data is
only prescribed in the half line $[a,\infty)$ and that it is given by
a piecewise constant function $u_o$ whose jumps are located in the
same positions as the jumps in $k$ and whose values satisfy
Rankine--Hugoniot conditions with zero speed.

With this particular problem, we are attempting to describe the case
of a physical system where some obstructions have appeared in the past
and then the evolution has stabilized into a stationary
solution. Using again the traffic flow model in
Example~\ref{ex:highway}, consider the case when
 an accident occurred in
the unobservable interval $(a,b)$ at some  time in the past. The accident 
caused all cars to slow down until they overtook the section of the
road obstructed by the vehicles involved, causing an increase in the
density of cars localized only in some interval
$J=[\xi_1,\xi_2]\subset\,(a,b)$, whose endpoints cannot be deduced
from the density of cars in $[b,\infty)$, where the accident does not
effect the dynamics. In this case, the only way to gather additional
information is to change the number of cars entering at $x=a$ and to
observe how this change affects the solution in the observable region
$[b,\infty)$. In other words, this problem could be considered as an
initial-boundary value problem in which we are free to choose suitable
boundary data $u_{\rm bdry}$ at $x=a-\ve$ for a fixed $\ve>0$, so that
the observations in $[0,T]\times\big([a-\ve,a]\cup[b,\infty)\big)$ of
the solution to the initial-boundary value problem
\begin{align*}
  &\partial_t u + \partial_x \big(k(x)f(u)\big) =0\,,
  \mbox{ in }[0,T]\times [a-\ve,\infty), \\
  &u(0,x)=u_o(x), \quad x\in[a,\infty),\qquad\qquad u(t,a-\ve)=u_{\rm
    bdry}(t), \quad t\in[0,T]
\end{align*}
allow the computation of $k(x)$.

However, such a problem can be reformulated in terms of an auxiliary
Cauchy problem in the whole $[0,T]\times\R$, in which we are allowed
to choose the initial data $u_o$ in $(-\infty,a)$ instead of $u_{\rm
  bdry}$. In this way, we are going to use the observed solution to
the Cauchy problem for~\eqref{eq:intrononh} with initial data
$$
u(0,x)=u_o(x),\quad x\in\R\,,
$$
in order to reconstruct $k(x)$.  A posteriori, if we denote by $\hat
u(t,x)$ the solution to such a Cauchy problem, $\hat
u(t,x)|_{[a-\ve,\infty)}$ provides a solution to the initial-boundary
value problem with boundary data\footnote{Here and in the following we
  use the convention that
  $\phi(a\pm)=\lim_{\epsilon\downarrow0}\phi(a\pm\epsilon)$.}  $u_{\rm
  bdry}(t)=\hat u(t,a-\ve+)$.

Our result for this first problem is that if the unobservable region
$[a,b]$ is nowhere fully congested and if the observation interval
$[0,T]$  is large enough, then we can choose a suitable
initial data in $(-\infty,a)$ to reconstruct uniquely the function
$k(x)$ in $[a,b]$, and hence in the whole $\R$ thanks to {\bf (H2)}.
\begin{theorem}\label{thm:3} Assume that the conservation law
  satisfies {\bf (H1)}--{\bf (H3)}, that $f$ is a known function, that
  the initial data $u_o(\dott)$ is a stationary solution on
  $[a,\infty)$ attaining values in $[u_1,u_2]$, and that for each
  choice of a $\BV$ initial data $u_o(\dott)$ in $(-\infty,a)$ with
  values in $[u_1,u_2]$, the solution $u_{\rm obs}(t,x)$ to the
  corresponding Cauchy problem for~\eqref{eq:intrononh} is partially
  observable, in the sense of Definition~\ref{def:observable}{\it (b)}.

  \n Then, if $u_o(\dott)<u_2$ in $[a,b]$, there exists $T>0$ large
  enough and a unique choice of $(k_1,\xi_1,\xi_2)$ such that, by
  denoting by $u_{(k_1,\xi_1,\xi_2)}$ the solution
  to~\eqref{eq:intrononh} with initial data $u_o$ and with $k(x)$
  defined in~\eqref{eq:disck}, there holds
  \begin{equation}\label{eq:reconstr_3}
    u_{(k_1,\xi_1,\xi_2)}(t,x)=u_{\rm obs}(t,x), \quad (t,x)\in[0,T]\times(\,\R\setminus\,(a,b)\,)\,. 
  \end{equation}
\end{theorem}

The actual reconstruction procedure for the piecewise constant function
$k(x)$ will be given in Section~\ref{sec:proof} (see in particular the proof 
of Lemma~\ref{lem:3-2} and Remark~\ref{rem:faster_proc}). Here, we want to comment about the assumption 
$u_o(\dott)<u_2$ in $[a,b]$.

\begin{remark}\label{rem:stationary_data} In Theorem~\ref{thm:3} we
  need the hypothesis that $u_o(x)< u_2$ for all $x\in[a,b]$, i.e.,
  that no part of the unobservable region is fully congested, to
  complete the reconstruction procedure. This assumption needs some
  comments in view of possible applications, because it appears to
  require information on the initial state of the physical system that
  cannot be known based only on partial observability.

  As a preliminary fact, note that the assumption that $u_o(\dott)$ is
  a stationary solution to~\eqref{eq:intrononh} in $[a,\infty)$,
  together with {\bf (H1)}--{\bf (H3)}, implies that
  \begin{equation}\label{eq:stationary_data}
    u_o(x)=\begin{cases}
      u_o(a),  & \mbox{if } a<x<\xi_1\,,\\
      \omega,  & \mbox{if } \xi_1<x<\xi_2\,,\\
      u_o(b),  & \mbox{if } x>\xi_2\,,
    \end{cases}
  \end{equation}
  for some constant $\omega\in [u_1,u_2]$, and that the jumps at
  $x=\xi_1$ and $x=\xi_2$ must be stationary. The Rankine--Hugoniot
  condition implies that there also holds
  \begin{equation}\label{eq:boundary_cond}
    k_o f(u_o(a))=k_1 f(\omega)=k_o f(u_o(b))\,,
  \end{equation}
  where the quantities $k_o,u_o(a),u_o(b)$ are known and the
  quantities $k_1,\omega$ are unknown.

  If $\{u_o(a),u_o(b)\}\subseteq \,(u_1,u_2)$, then $f>0$ in all the
  above equalities~\eqref{eq:boundary_cond} and also $\omega<u_2$ must
  hold. As a result, no part of the region $(a,b)$ can be fully
  congested and Theorem~\ref{thm:3} can be applied.

  In the case of either $u_o(a)=u_2$ or $u_o(b)=u_2$, the
  reconstruction procedure cannot be applied; indeed, it would be
  impossible to reconstruct the value $k_1$ attained by $k(x)$ in the
  interval $[\xi_1,\xi_2]$, because~\eqref{eq:boundary_cond} simply
  implies $f(\omega)=0$ independently of $k_1$. On a positive note,
  however, such an impossibility can also be immediately detected by
  the known values of $u_o(\dott)$ in $x=a$ or $x=b$.

  It remains to consider the case of $u_o(a)=u_o(b)=u_1$. In this
  case, Theorem~\ref{thm:3} applies if $\omega=u_1$ and fails if
  $\omega=u_2$. Since we cannot observe $u_o(\dott)$ in
  $[\xi_1,\xi_2]$, it is not a priori possible to decide in which case
  we are.  Trying to apply the reconstruction procedure to a problem
  where $u_o(x)=u_2$ in $[\xi_1,\xi_2]$ soon leads to the appearance
  of the ``forbidden'' state $u=u_2$ at $x=a$ so that the
  assumption of partial observability allows us   a posteriori  to detect
  the presence of a fully congested region inside $(a,b)$.

  Therefore, in the case $u_o(a)=u_o(b)=u_1$, which is the only one in
  applications where it would be impossible to know in advance if the
  assumption $u_o(\dott)<u_2$ is satisfied, the conclusion of the
  theorem could be reformulated as follows: either there exists
  $T_1>0$ such that $u_{\rm obs}(T_1,a+)=u_2$, or there exists $T_2>0$
  large enough and a unique choice of $(k_1,\xi_1,\xi_2)$ such that
  \begin{equation*}
    u_{(k_1,\xi_1,\xi_2)}(t,x)=u_{\rm obs}(t,x),\quad (t,x)\in[0,T_2]\times(\R\setminus\,(a,b))\,. 
  \end{equation*}
\end{remark}
\begin{remark} In the application to the traffic flow model in
  Example~\ref{ex:highway}, Theorem~\ref{thm:3} covers for instance
  the case of initial data $u_o(\dott)\equiv 0$ in $[a,\infty)$. In
  other words, among other cases, the reconstruction procedure in
  Theorem~\ref{thm:3} allows us to recover $k(x)$ when we are
  considering a highway which is known to be empty at $t=0$.
\end{remark}

The second problem we consider, under the assumption of partial
observability, is the reconstruction of $k(x)$ in the case
of initial data $u(0,x)\equiv \bar u_o$, for some constant $\bar
u_o\in\,[u_1,u^m)$. In other words, we assume that at the initial time
the whole spatial domain contains a constant state $\bar u_o$. With
this particular problem, we are attempting to describe the case of a
physical system in which, at time $t=0$, the constant flux function
$k_{\rm old}(x)\equiv k_o$ is suddenly replaced by a piecewise
constant function $k(x)$, due to the appearance of some obstructions
in the system. Considering once again the traffic flow model in
Example~\ref{ex:highway}, you might think of a constant density of
cars distributed in the whole highway and of a car accident occurring,
at time $t=0$, in some place inside the unobservable interval $(a,b)$.

In this case, the initial data $u\equiv \bar u_o$ is not a stationary
solution for~\eqref{eq:intrononh} with discontinuous flux $k(x)f(u)$
and therefore the solution will immediately develop additional waves
around the discontinuity points for $k$.

Our result for this problem is that, if we observe the solution long
enough, then we can always reconstruct the function $k(x)$ in $[a,
b]$, and hence in the whole $\R$, as before (details on the reconstruction are
given in Section~\ref{sec:proof}, in particular in the proof of Lemma~\ref{lem:4-2}). 
Uniqueness of the resulting flux $k(x)$, on the other hand, only holds when the
obstruction is large enough. This is not entirely surprising, because
it is expected that the effect of a very small obstruction occurring
in a very small spatial region $[\xi_1,\xi_2]\subseteq[a,b]$ gets
canceled before reaching the observable region
$\R\setminus\,(a,b)$. But it might also happen that the obstruction
produces effects that can be detected in the observable region and
still the data is insufficient to lead to a unique reconstruction: in
the latter case, it is in general possible to provide infinitely many
functions $k(x)$, all leading to the same solution in
$\R\setminus\,(a,b)$.

\begin{theorem}\label{thm:4} Assume that the conservation law satisfies
  {\bf (H1)}--{\bf (H3)}, that $f$ is a known function, and that the
  solution $u_{\rm obs}(t,x)$ to the Cauchy problem
  for~\eqref{eq:intrononh} with a constant initial data
  $u(0,\dott)\equiv \bar u_o\in\,[u_1,u^m)$ is partially observable, 
  in the sense of Definition~\ref{def:observable}{\it (a)}.

  Then, either $u_{\rm obs}(t,x)\equiv \bar u_o$ for all
  $(t,x)\in[0,\infty)\times(\R\setminus\,(a,b))$, and hence we can
  assume $k(x)\equiv k_o$ for all $x\in\R$, or there exist $T > 0$
  large enough and a choice of $(k_1,\xi_1,\xi_2)$ such that, denoting
  $u_{(k_1,\xi_1,\xi_2)}$ the solution to~\eqref{eq:intrononh} with
  $k(x)$ defined in~\eqref{eq:disck}, there holds
  \begin{equation}\label{eq:reconstr_4}
    u_{(k_1,\xi_1,\xi_2)}(t,x)=u_{\rm obs}(t,x),\quad (t,x)\in[0,T]\times(\R\setminus\,(a,b))\,.
  \end{equation}
  Moreover, if there exists $T_1\in\,(0,T)$ such that
  \begin{equation*}
    u(T_1,a-)=\bar u_o< u(T_1,a+)\,,
  \end{equation*}
  or
  \begin{equation*}
    u(T_1,b+)=\bar u_o> u(T_1, b-)~~~\mbox{ and
    }~~~\inf\big\{s\in\,(T_1,T)~;~u(s,b)>u(T_1,b-)\big\}> T_1\,,
  \end{equation*}
  then the choice is unique.
\end{theorem}
A few comments are in order. First of all, we notice that
Theorems~\ref{thm:3} and~\ref{thm:4} state that there exists an
observation time $T>0$ large enough so that the reconstruction
procedure can be completed successfully. The reason for this is that
we need enough waves to pass through the unobservable region $[a,b]$
and reach the observable region, before we can fully determine
$k$. If, e.g., the constant $k_1$ in~\eqref{eq:disck} is close to
zero, the waves can take a very long time $\overline{T} \approx
{\mathcal O}(1) {b-a\over k_1}$ to pass through the unobservable
region and therefore the reconstruction is not possible by only
observing the solution in $[0,\tau]$ with $\tau<\overline{T}$. The
technical Lemmas~\ref{lem:3-3} (for Theorem~\ref{thm:3})
and~\ref{lem:4-2} (for Theorem~\ref{thm:4}) show the properties
satisfied by the observed solution $u_{\rm obs}$ at time
$\overline{T}$, and characterize the minimal time $T$ for which the
reconstruction procedure can be completed.

Also, we want to emphasize some features of our results. % which should
%be taken into account in applications to practical problems. 
In Theorem~\ref{thm:1} several solutions corresponding to Riemann initial
data have to be observed, but only a single observation for each solution 
(at time $T>0$) is needed. Since $T$ can be chosen arbitrarily small, we can
always test a large number of initial data for a very short time, so to
obtain an accurate  piecewise affine approximation of the flux $f(u)$ 
in time smaller than any fixed $\hat T>0$.
In Theorem~\ref{thm:2} observations have to be
performed on a whole interval $(0,T]$ to recover the flux function
$k(x)$. Once again, however, $T$ can be chosen arbitrarily small and 
an exact reconstruction of $k(x)$ can be found in arbitrarily small time.

In both cases above, the reconstruction in small time is possible 
because we can observe the solution on the whole spatial
domain and because we are free to select any initial data. 
Although this procedure might not be immediate to apply in 
practical situations because, for instance, setting up multiple 
initial data requires time and efforts, it may still be of help in cases
of bounded propagation speed, when observations on a small 
spatial intervals are enough to reconstruct the flux: in these situations
one could observe different initial data in different portions of the road,
so to actually improve the reconstruction while reducing the 
necessary efforts.

On the other hand, when the solution cannot be observed in the whole $\R$, as in
Theorems~\ref{thm:3} and~\ref{thm:4}, it becomes vital to study
$u_{\rm obs}$ on an interval $[0,T]$, with $T$ possibly
very large as remarked earlier. At the same time, the choice of initial
data becomes more important, because carefully chosen initial data can
convey more information about the flux.

It is not surprising, therefore, that when there are unobservable
regions and we can observe the solution corresponding to any
initial data of our choice, as in Theorem~\ref{thm:3}, we still can
recover a unique exact reconstruction of $k(x)$, unless the
unobservable region is fully congested.
But in cases when there are unobservable regions and the initial
data cannot be freely chosen, like in Theorem~\ref{thm:4}, the amount
of information that can be recovered from the solution is limited. In
particular, in some cases we lose the uniqueness of the reconstructed
flux $k(x)$, because the effect on the given initial data of many
different small obstructions might pass equally undetected in the
observable region.

Finally, we remark that the assumption {\bf (H3)} on $k(x)$, by
prescribing the exact number of discontinuities in the unobservable
region, is very strong. However, {\bf (H3)} is really necessary for
the inverse problem to be well-posed: in Section~\ref{sec:counterex}
we present a few examples where, one by allowing for three or more
jumps in $k(x)$, immediately is led to the existence of infinitely
many piecewise constant functions $\hat k(x)$, whose corresponding
solutions coincide with $u_{\rm obs}$ in $[0,T]\times(\R\setminus
[a,b])$. In other words, the reconstruction problem is in general
ill-posed within the class of piecewise constant functions which do
not satisfy {\bf (H3)}.

\section{Ill-posedness when $k(x)$ has more than two
  jumps}\label{sec:counterex} 
In this section, we show through a few examples that the problem with
partial observability is in general ill--posed whenever the
function $k(x)$ in~\eqref{eq:intrononh} is allowed to have three or
more jumps, i.e.,  when $k$ satisfies {\bf (H1)} and {\bf (H2)} but not
{\bf (H3)}. Namely, we show that in several situations there exist infinitely many
different functions $k$ with three or more jumps which produce exactly the same 
solution in the observable region $\R\setminus\,(a,b)$. Considering 
the car traffic example, this means that in some situations there
could be $2$ or $3$ or more small accidents in the region $(a,b)$ or a single larger one, 
and there would be no way to distinguish between them by just observing the 
situation in $\R\setminus\,(a,b)$. 
This is always the case, for instance, if the accident which is closer to the extreme 
$x=a$ reduces the flux more than the subsequent ones.

\n In view of these examples, and of the fact that there is no 
reason in applications to exclude obstructions which are larger close to
$x=a$ than in the rest of the region, one can think to hypothesis {\bf (H3)} as
a way to single out an approximation of the real, and possibly very complex, 
structure of $k(\dott)$ in $(a,b)$ by means of a single obstruction. In turn, this 
approximation is ``good'' because the corresponding solution in $\R\setminus (a,b)$
coincides with the observed one for all times if the flux function satisfies {\bf (H3)}
or if we are in any of the cases below.

\begin{example}\label{ex:counter1} Consider the Cauchy problem
  for~\eqref{eq:intrononh} with initial data $u(0,x)\equiv u_1$ for
  $x\in \,(a,\infty)$. In the highway Example~\ref{ex:highway}, this
  initial data means that the road is initially empty. Assume that we
  have reconstructed a coefficient $\kappa(x)$ so that the solution
  $u^\kappa$ to $\partial_t u +\partial_x(\kappa(x)f(u))=0$ coincides
  with $u_{\rm obs}$ in the observable region
  $[0,T]\times\big(\R\setminus\,(a,b)\big)$, and that
  \begin{equation}\label{eq:flux_candidate}
    \kappa(x) =
    \left\{
      \begin{array}{ll}
        k_o, & x \notin \,(\xi', \xi' + \chi_1 + \chi_2] \,, \\
        k_1, & x\in \,(\xi', \xi' + \chi_1] \,, \\
        k_2, & x\in \,(\xi' + \chi_1, \xi' + \chi_1 + \chi_2] \,,
      \end{array}\right.
  \end{equation}
  for suitable positive numbers $\chi_1, \chi_2$ such that
  $\chi_1+\chi_2 \le b - a$, for $0<k_1< k_2<k_o$ and for a fixed
  $\xi' \in [a, b - \chi_1 - \chi_2]$.

  If $\chi_1+\chi_2 < b - a$, we claim that for every $\ve>0$ small
  enough, also the solutions $u^{\kappa_\ve}$ coincide with $u_{\rm
    obs}$ in $[0,T]\times\big(\R\setminus\,(a,b)\big)$ if we choose
  the coefficient $\kappa_\ve$ as follows (see Figure~\ref{fig:counterex1}, middle)

  \begin{figure}\begin{center}
      \includegraphics[width=.99\textwidth]{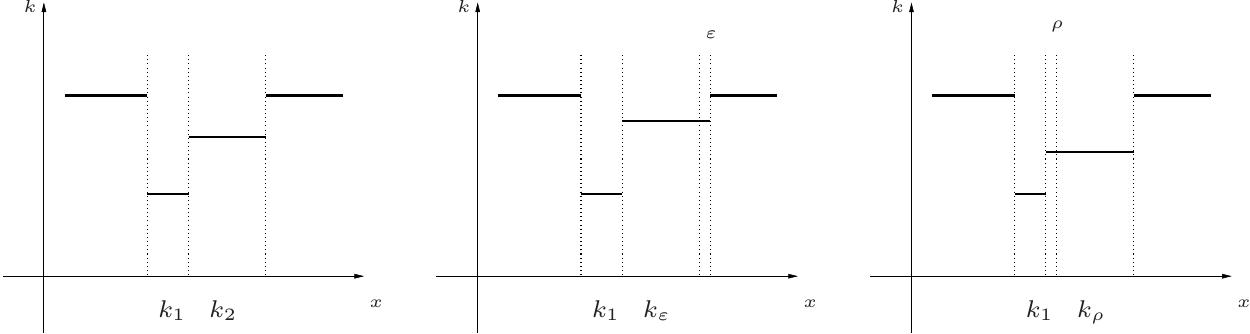}
      \caption{Three different $k(x)$ producing the same solution
        outside $(a,b)$, as in Examples~\ref{ex:counter1} and~\ref{ex:counter2}.\label{fig:counterex1}}
  \end{center}\end{figure}

  \begin{equation*}
    \kappa_\ve(x) = \begin{cases}
      k_o, & x \notin \,(\xi', \xi' + \chi_1 + \chi_2+\ve] \,, \\
      k_1, & x\in \,(\xi', \xi' + \chi_1] \,, \\
      k_\ve, & x\in \,(\xi' + \chi_1, \xi' + \chi_1 + \chi_2+\ve] \,,
    \end{cases}
  \end{equation*}
  with $k_\ve\in\,(k_2,k_o)$ given by
  \begin{equation*}
    k_\ve=\,{\chi_2+\ve\over\displaystyle{\chi_2\over k_2}\,+\,{\ve\over
        k_o}}\, = k_2\left(1+\ve\,\,{\displaystyle{1\over k_2}\,-\,{1\over
          k_o}\over\displaystyle{\chi_2\over k_2}\,+\,{\ve\over
          k_o}}\right)\,.
  \end{equation*}
  Indeed, independently of the choice of the initial data $u(0,\dott)$
  in $(-\infty,a]$ (or of the boundary data $u_{\rm bdry}$ at $x=a$)
  the exact same solution will always be observed for $x\geq b$. This
  can be seen as follows. Fix any initial data in $(-\infty,a]$
  attaining some value $\omega$ larger than $u_1$. Then the solution
  at time $t=0+$ will contain a centered rarefaction wave traveling
  with speed $k_o f'(u)$ for $u\in [u_1,\omega]$. In particular, by
  {\bf (H1)} states close to $u_1$ will travel with positive speed,
  i.e., towards the unobservable region $[a,b]$. Eventually, the
  centered rarefaction wave will cross completely $[a,b]$ and emerge
  at $x=b$ after having spent in $[a,b]$ a time
  \begin{equation*}
    T_{u}={b - a-(\chi_1+\chi_2) \over k_o f'(u)} + {\chi_1 \over
      k_1f'(u)} + {\chi_2 \over k_2 f'(u)}\,,\qquad\qquad u>u_1\,,
  \end{equation*}
  if the flux is $\kappa f$, and a time
  \begin{equation*}
    T'_{u}={b - a-(\chi_1+\chi_2+\ve) \over k_o f'(u)} + {\chi_1 \over
      k_1f'(u)} + {\chi_2+\ve \over k_\ve f'(u)}\,,\qquad\qquad u>u_1\,,
  \end{equation*}
  if the flux is $\kappa_\ve f$. It is easy to verify that the choice
  of $k_\ve$ implies $T_u=T'_u$ for all states $u$ which pass $[a,b]$,
  proving that the solution restricted to $\R\setminus\,(a,b)$ is the
  same for both fluxes. Therefore, any function $\kappa_\ve(x)$
  provides a solution to our inverse problem.
\end{example}

\begin{example}\label{ex:counter2} We now show that the loss of uniqueness
  cannot be avoided by prescribing the length of the ``obstruction''
  interval $\chi_1+\chi_2$. Indeed, let us consider the same problem
  as in Example~\ref{ex:counter1} and the same possible flux function
  $\kappa(x)$ defined in~\eqref{eq:flux_candidate} for suitable
  positive numbers $\chi_1, \chi_2$ such that $\chi_1+\chi_2 \le b -
  a$, and for $0<k_1< k_2<k_o$. It can be easily verified that, for
  any fixed $\rho>0$ small enough, the solutions corresponding to flux
  functions (see Figure~\ref{fig:counterex1}, right)

  \begin{equation*}
    \kappa_\rho(x) = \begin{cases}
      k_o, & x \notin \,(\xi', \xi' + \chi_1 + \chi_2] \,, \\
      k_1, & x\in \,(\xi', \xi' + \chi_1-\rho] \,, \\
      k_\rho, & x\in \,(\xi' + \chi_1-\rho, \xi' + \chi_1 + \chi_2] \,,
    \end{cases}
  \end{equation*}
  with $k_\rho\in\,(k_1,k_2)$ defined by
  \begin{equation*}
    k_\rho=\,{\chi_2+\rho\over\displaystyle{\chi_2\over
        k_2}\,+\,{\rho\over k_1}}\, =
    k_2\left(1-\rho\,\,{\displaystyle{1\over k_1}\,-\,{1\over
          k_2}\over\displaystyle{\chi_2\over k_2}\,+\,{\rho\over
          k_1}}\right)
  \end{equation*}
  once again coincide on $\R\setminus\,(a,b)$ with the ones found in
  Example~\ref{ex:counter1}.
\end{example}

\begin{example}\label{ex:counter3} The previous examples can be easily
  generalized to the case of a flux function with four or
  more discontinuities. For instance, assuming that we have
  reconstructed the function
  \begin{equation*}
    \kappa(x) = \begin{cases}
      k_o, & x \notin \,(\xi', \xi' + \chi_1 + \chi_2+\chi_3) \,, \\
      k_1, & x\in \,(\xi', \xi' + \chi_1) \,, \\
      k_2, & x\in \,(\xi' + \chi_1, \xi' + \chi_1 + \chi_2) \,,\\
      k_3, & x\in \,(\xi' + \chi_1 + \chi_2, \xi'  + \chi_1 + \chi_2+\chi_3) \,,
    \end{cases}
  \end{equation*}
  for suitable positive constants $\chi_1, \chi_2, \chi_3$ such that
  $\sum_i \chi_i \le b - a$ and for $0 < k_1 < k_2 < k_3 < k_o$, one
  can easily prove that the solutions to~\eqref{eq:intrononh} with
  flux $\kappa f$ and initial data $u(0,x)\equiv u_1$ for $x\in
  \,(a,\infty)$ coincide (outside $(a,b)$) with the solutions
  to~\eqref{eq:intrononh} with flux $\kappa_\ve f$ and the same
  initial data, if we define
  \begin{equation*}
    \kappa_\ve(x) = \begin{cases}
      k_o, & x \notin \,(\xi', \xi' +\sum_i \chi_i) \,, \\
      k_1, & x\in \,(\xi', \xi' + \chi_1) \,, \\
      \hat k_\ve, & x\in \,(\xi' + \chi_1, \xi' + \chi_1 + \chi_2-\ve) \,,\\
      \tilde k_\ve, & x\in \,(\xi' + \chi_1 + \chi_2-\ve, \xi'  + \sum_i \chi_i) \,,
    \end{cases}
  \end{equation*}
  for $\ve>0$ small enough and for any $\hat k_\ve,\tilde k_\ve$ such
  that
  \begin{equation*}
    {\chi_3\over k_3}\,+\,{\chi_2\over k_2} \,=\,{\chi_3+\ve\over \tilde
      k_\ve}\,+\,{\chi_2-\ve\over \hat k_\ve}\,.
  \end{equation*}
  In particular, by choosing
  \begin{equation*}
    \tilde k_\ve=\hat k_\ve=\ell:= k_2\left({\chi_2/k_2\over
        \chi_2/k_2 +\chi_3/k_3}\right)+k_3\left({\chi_3/k_3\over
        \chi_2/k_2 +\chi_3/k_3}\right)\in\,(k_2,k_3)
  \end{equation*}
  one obtains that the same solution $u_{\rm obs}$ in
  $\R\setminus\,(a,b)$ corresponding to $\kappa(x)f(u)$ can also be
  obtained as solution of the conservation law with flux $\bar
  \kappa(x)f(u)$ where
  \begin{equation*}
    \bar \kappa(x) = \begin{cases}
      k_o, & x \notin \,(\xi', \xi' +  \sum_i \chi_i) \,, \\
      k_1, & x\in \,(\xi', \xi' + \chi_1) \,, \\
      \ell, & x\in \,(\xi' + \chi_1, \xi' +  \sum_i \chi_i) \,,
    \end{cases}
  \end{equation*}
  i.e., not only the available data are insufficient to distinguish
  between flux functions with a different number of discontinuities
  (in this case three or four jumps), but it is possible to construct
  infinitely many additional flux functions by applying the ideas of
  Examples~\ref{ex:counter1} and~\ref{ex:counter2} to $\bar\kappa$,
  and all these fluxes would give solutions coinciding with $u_{\rm
    obs}$ in the observable region $\R\setminus\,(a,b)$.
    
  \n Finally, we remark that by repeating the same argument on $\bar \kappa$ and by defining
  \begin{equation*}
    \ell':= k_1\left({\chi_1/k_1\over
        \chi_1/k_1 +(\chi_2+\chi_3)/\ell}\right)+\ell\left({(\chi_2+\chi_3)/\ell\over
        \chi_1/k_1 +(\chi_2+\chi_3)/\ell}\right)\in\,(k_1,\ell)
  \end{equation*}
  and
  \begin{equation*}
    \bar \kappa'(x) = \begin{cases}
      k_o, & x \notin \,(\xi', \xi' +  \sum_i \chi_i) \,, \\
      \ell', & x\in \,(\xi', \xi' +  \sum_i \chi_i) \,,
    \end{cases}
  \end{equation*}
  one also obtains a flux function which satisfies {\bf (H3)} and produces the
  same solution $u_{\rm obs}$ in the observable region $\R\setminus\,(a,b)$.
\end{example}

\begin{example}\label{ex:counter4} As a last example of
  ill-posedness, we show that when four or more discontinuities are
  assumed to be present in $k(x)$, then not even imposing a priori the
  length of each discontinuity helps to recover uniqueness.

  Once again, consider the Cauchy problem for~\eqref{eq:intrononh}
  with the initial data $u(0,x)\equiv u_1$ for $x\in \,(a,\infty)$.
  Fix three positive numbers $\chi_1, \chi_2, \chi_3$ such that
  $\sum_i \chi_i \le b - a$, representing the length of the intervals
  in which $k(x)\neq k_o\equiv k|_{\R\setminus\,(a,b)}$ as in
  Example~\ref{ex:counter3}, and fix $\xi' \in [a, b - \sum_i \chi_i]$
  representing the location of the first discontinuity of $k$. Assume
  that we reconstruct a flux function (see Figure~\ref{fig:counterex2})
  \begin{equation*}
    \kappa_1(x) = \begin{cases}
      k_o, & x \notin \,(\xi', \xi' + \sum_i \chi_i) \,, \\
      k_1, & x\in \,(\xi', \xi' + \chi_1) \,, \\
      k_2, & x\in \,(\xi' + \chi_1, \xi' + \chi_1 + \chi_2) \,,\\
      k_3, & x\in \,(\xi' + \chi_1 + \chi_2, \xi' + \sum_i \chi_i) \,,
    \end{cases}
  \end{equation*}
  with $0 < k_1 < k_2 < k_3 < k_o$, so that the solution
  $u^{\kappa_1}$ to $\partial_t u +\partial_x(\kappa_1(x)f(u))=0$
  coincides with $u_{\rm obs}$ in the observable region
  $[0,T]\times\big(\R\setminus\,(a,b)\big)$. Then it is easy to verify
  that also the piecewise constant function defined by
  % \begin{figure}\begin{center}
  %     \psfrag{x}{${\scriptscriptstyle x}$}
  %     \psfrag{k}{${\scriptscriptstyle k}$} \psfrag{1}{${\scriptstyle
  %         k_1}$} \psfrag{2}{${\scriptstyle k_2}$}
  %     \psfrag{3}{${\scriptstyle k_3}$} \psfrag{a}{${\scriptstyle
  %         k_1}$} \psfrag{b}{${\scriptstyle k_3}$}
  %     \psfrag{c}{${\scriptstyle k_2}$}
  %     \includegraphics[width=.8\textwidth]{multi.eps}\\ Figure~2:
  %     Two coefficients $k(x)$ producing the same solution outside
  %     $(a,b)$.
  %   \end{center}\end{figure}
  \begin{figure}\begin{center}
      \includegraphics[width=.9\textwidth]{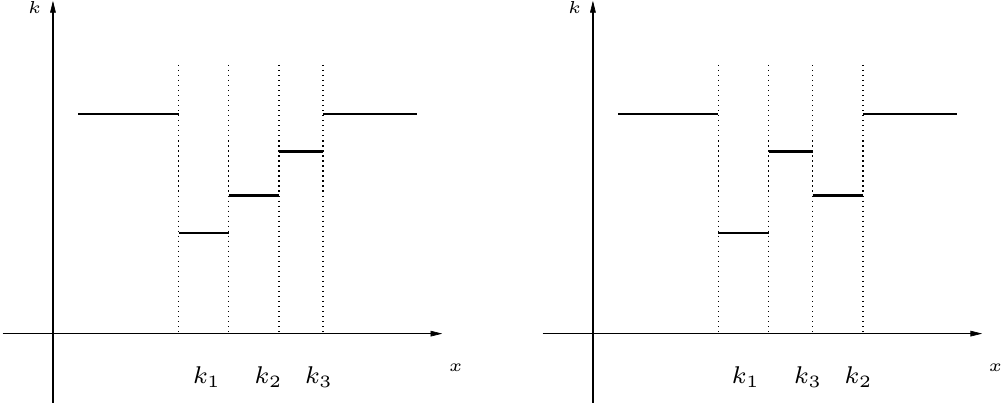}
      \caption{Two choices of $k(x)$ with four jumps which produce the same solution
        outside $(a,b)$, as in Example~\ref{ex:counter4}.\label{fig:counterex2}}
    \end{center}\end{figure}
  \begin{equation*}
    \kappa_2(x) = \begin{cases}
      k_o, & x \notin \,(\xi', \xi' + \sum_i \chi_i) \,, \\
      k_1, & x\in \,(\xi', \xi' + \chi_1) \,, \\
      k_3, & x\in \,(\xi' + \chi_1, \xi' + \chi_1 + \chi_3) \,,\\
      k_2, & x\in \,(\xi' + \chi_1 + \chi_3, \xi' + \sum_i \chi_i) \,,
    \end{cases}
  \end{equation*}
  i.e., obtained by switching the interval where $\kappa_1(x)=k_2$ and
  $\kappa_1(x)=k_3$, gives a solution $u^{\kappa_2}$ which coincides
  with $u_{\rm obs}$ in $\R\setminus\,(a,b)$.
\end{example}

\section{Technical proofs and reconstruction procedures}\label{sec:proof}
\subsection*{Proof of Theorem~\ref{thm:1}}
As already remarked in Section~\ref{sec:main_res}, once we have
proved~\eqref{eq:recest}, the general stability result \cite[Theorem
2.13]{HoldenRisebroBook} ensures that also~\eqref{eq:close} is
satisfied. Hence, the proof reduces to the construction of the
approximated flux $f_\nu$ which satisfies~\eqref{eq:recest}.

Let $T>0$, $u_*,u^*\in\R$ such that $u_*<u^*$ and $c\in\R$ be fixed.
Fix also $\nu\in\N$ and define $\delta$ and $\{u_o,\ldots,
u_{2^\nu}\}$ as in the statement of Theorem~\ref{thm:1}. Of course, we start by
definining $f_\nu(u_*)=c$.

Since we only assume to be able to observe the solution
to~\eqref{eq:introhom} at time $T$, we have to choose carefully the
initial data. In this case, let us consider the following family of
Riemann data:
\begin{equation}\label{eq:riemdata}
  u_o^h(x)=\left\{
    \begin{array}{ll}
      u_h, & \,\,x<0,\\
      u_{h+1}, & \,\,x>0,
    \end{array}
  \right.
  \qquad h=0,\ldots,2^\nu-1.
\end{equation}
The strategy is to use the solution corresponding to each $u_o^h$ to
assign $f_\nu$ in $u_{h+1}$. First, we consider the particular case of
a solution which at time $T$ consists of a single wave, either an
entropy shock wave or a centered rarefaction wave, joining the states
$u_h$ and $u_{h+1}$. This is the case, for instance, when $f$ has no inflection
points in the interval $[u_h, u_{h+1}]$. Once we know how to deal with
this easier case, we move to the general situation.

\n {\bf Step 1 (Shock).}  Fixed $h\ge 0$, let $\tilde u(\dott)=u_{\rm
  obs}(T,\dott)$ be the solution
to~\eqref{eq:introhom}--\eqref{eq:riemdata} at time $T$, consisting of
a single shock wave joining $u_h$ and $u_{h+1}$, and let $x_h\in\R$ be
the location of the jump. Then, the propagation speed of this wave is
given by $s_h=x_h/ T$ and, by Rankine--Hugoniot conditions, there
holds
\begin{equation*}
  f(u_{h+1})=f(u_h)+s_h(u_{h+1}-u_h)=f(u_h)+\,{\delta x_h\over T}\,.
\end{equation*}
Therefore, if $f_\nu(u_o),\ldots, f_\nu(u_h)$ are given so that
$f_\nu(u_\alpha)=f(u_\alpha)$, we can define
\begin{equation*}
  f_\nu(u_{h+1})=f_\nu(u_h)+\,{\delta x_h\over T}\,=f(u_h)+\,{\delta
    x_h\over T}\,=f(u_{h+1})\,.
\end{equation*}

\n {\bf Step 2 (Rarefaction).}  Fixed $h\ge 0$, let $\tilde
u(\dott)=u_{\rm obs}(T,\dott)$ be the solution
to~\eqref{eq:introhom}--\eqref{eq:riemdata} at time $T$, consisting of
a single centered rarefaction wave joining $u_h$ and $u_{h+1}$, and
let $I_h=[x_h, x_{h+1}]$ be the interval in which $\tilde u(\dott)$ is
not constant. Then, if $f_\nu(u_o),\ldots, f_\nu(u_h)$ are given so
that $f_\nu(u_\alpha)=f(u_\alpha)$, we claim that
\begin{itemize}
\item if we replace this rarefaction wave with a shock wave separating
  the same states, whose jump is located at the point
  \begin{equation}\label{eq:rarjump}
    \xi_h:= {\int_{u_h}^{u_{h+1}} x(u)\,du \over u_{h+1} -
      u_h}={\int_{u_h}^{u_{h+1}} x(u)\,du \over \delta} 
  \end{equation}
  where $u\mapsto x(u)$ is the inverse of $x\mapsto \tilde u(x)$ on
  $I_h$, and traveling with speed $\xi_h /T$;
\item and, if we define
  \begin{equation}\label{eq:rarflux}
    f_\nu(u_{h+1})=f_\nu(u_h)+\,{\delta \xi_h\over T},
  \end{equation}
\end{itemize}
then $f_\nu(u_{h+1})=f(u_{h+1})$.  Indeed, it is enough to recall
that, by definition of centered rarefaction waves (see, e.g.,
\cite{BressanBook}), the following equality holds
\begin{equation*}
  \int_{u_h}^{u_{h+1}} {x(u) \over T}\,du = \int_{u_h}^{u_{h+1}}
  f'(u)\,du\,.
\end{equation*}
Therefore, we have
\begin{align*}
  f_\nu(u_{h+1})&=f_\nu(u_h)+\,{\delta \xi_h\over T}=
  f(u_h)+\int_{u_h}^{u_{h+1}} {x(u) \over T}\,du\\ &=
  f(u_h)+\int_{u_h}^{u_{h+1}} f'(u)\,du = f(u_{h+1})\,.
\end{align*}
Note that the computation of the point $\xi_h$ used to define
$f_\nu(u_{h+1})$ can be done explicitly under the observability
assumption. Indeed, $\tilde u(\dott)$ is monotonically
increasing in $I_h$. Therefore, once we know $\int_{x_h}^{x_{h+1}}
\tilde u(x)dx$ and the value attained by $\tilde u$ at the points
$x_h,x_{h+1}$, then we also know the value of the integral
$\int_{u_h}^{u_{h+1}} x(u)\,du$ used to define $\xi_h$ (see
Lemma~\ref{lem:inverse_int}), even without computing the expression of
the inverse function $u\mapsto x(u)$.

\n {\bf Step 3 (General case).}  Fixed $h\ge 0$, assume
$f_\nu(u_o),\ldots, f_\nu(u_h)$ are given so that
$f_\nu(u_\alpha)=f(u_\alpha)$ and let $\tilde u(\dott)=u_{\rm
  obs}(T,\dott)$ be the solution
to~\eqref{eq:introhom}--\eqref{eq:riemdata} observed at time $T$. In
general, $\tilde u$ can consist of more than one single wave but, in
any case, $\tilde u$ is monotonically increasing and it can contain
only a finite number of different waves, because of the choice of the
initial data and because we are assuming that $f$ has a finite number of 
inflection points.

Therefore, let $x_1<\dots<x_{M_1}$ be the locations of jumps of
$\tilde u$ and let $I_1,\ldots, I_{M_2}$ the intervals in which
$\tilde u$ has non-zero derivative. To define $f_\nu(u_{h+1})$, we
simply proceed applying the construction in Step~1 to each shock and
the one in Step~2 to each centered rarefaction wave which appears in
$\tilde u(\dott)$.

Namely, we first replace each rarefaction joining two states $u^\ell,
u^r$ on the interval $I_j$ ($j=1,\ldots,M_2$) with a shock centered at
\begin{equation*}
  \xi_j={\int_{u^\ell}^{u^r} x(u)\,du \over u^\ell - u^r} \in I_j\,.
\end{equation*}

In this way, we obtain a new piecewise constant function $\bar
u(\dott)$ whose jumps are located at points $y_1<\dots<y_M$, with
$M=M_1+M_2$ and $\{y_1,\ldots, y_M\}=\{x_1,\ldots,
x_{M_1},\xi_1,\ldots, \xi_{M_2}\}$. Let $v_1<\dots<v_{M+1}$ be the
values attained by $\bar u$, i.e., let us assume
\begin{equation*}
  \bar u(x) =\begin{cases}
    v_1=u_h, & \mbox{if } x<y_1\,, \\
    \ \vdots & \\
    v_\alpha, &  \mbox{if } y_{\alpha-1}<x<y_{\alpha}\,,~\alpha=2,\ldots,M\,,  \\
    \ \vdots & \\
    v_{M+1}=u_{h+1}, &  \mbox{if } x>y_M\,. 
  \end{cases}
\end{equation*}
By construction,
\begin{itemize}
\item $\bar u(\dott)$ coincides with $\tilde u(\dott)$ outside
  $\bigcup_{k=1}^{M_2} I_k$;
\item on each interval $I_k$ where $\tilde u$ has a rarefaction
  joining $u^\ell, u^r$, $\bar u$ attains only the values $u^\ell,
  u^r$ and it jumps from $u^\ell$ to $u^r$ at a point $\xi_k\in I_k$
  such that
  \begin{equation*}
    {f(u^\ell)-f(u^r) \over u^\ell - u^r}\, = \,
    {\xi_k \over T}
  \end{equation*}
  as follows from~\eqref{eq:rarjump}--\eqref{eq:rarflux};
\item all the values $v_1,\ldots, v_{M+1}$ are known, since they are
  attained by $\tilde u(\dott)$ as adjacent states to shocks and
  rarefactions.
\end{itemize}
Now, set
\begin{equation*}
  y:= \sum_{\alpha=1}^M {v_{\alpha+1} - v_\alpha \over
    v_{M+1}-v_1}\,y_\alpha= \sum_{\alpha=1}^M {v_{\alpha+1} - v_\alpha
    \over \delta}\,y_\alpha
\end{equation*}
and define
\begin{equation*}
  f_\nu(u_{h+1})= f_\nu(u_h)+\,{\delta y \over T}\,.
\end{equation*}
We claim that $f_\nu(u_{h+1})=f(u_{h+1})$. Indeed,
\begin{align*}
  f_\nu(u_{h+1})&= f_\nu(u_h)+\,{\delta y \over T}=
  f(u_h)+\sum_{\alpha=1}^M {v_{\alpha+1} - v_\alpha \over T}\,
  y_\alpha\\ &= f(u_h)+\sum_{\alpha=1}^M[ f(v_{\alpha+1}) -
  f(v_\alpha)]\\ &= f(u_h)+f(v_{M+1}) - f(v_1)=f(u_{h+1})\,,
\end{align*}
where we have again used the Rankine--Hugoniot conditions and the
particular choices of $\xi_1,\ldots,\xi_{M_2}$ as locations for the
jumps in $\bar u$, which replace rarefactions in $\tilde u$.

At this point, we define $f_\nu$ on $[u_*,u^*]$ as the piecewise
affine function joining the values obtained in the previous steps:
\begin{equation*}
  f_\nu(u):= f_\nu(u_h)+ {f_\nu(u_{h+1}) - f_\nu(u_h) \over
    \delta}\,(u-u_h),\qquad u\in[u_h,u_{h+1}]\,.
\end{equation*}
Finally, we are ready to prove~\eqref{eq:recest}. Given any point
$u\in [u_*,u^*]$, there exists $\alpha\in\{0,\ldots, 2^\nu-1\}$ such
that $u\in~[u_\alpha, u_{\alpha+1}]$. Setting
$v_o=u_\alpha<\dots<v_N=u_{\alpha+1}$ the values such that $f$ is of
class $C^{1,1}$ on each interval $(v_j,v_{j+1})$, we then have
\begin{align*}
  |f'_\nu(u)-f'(u)|&=\Bigl| \,{f(u_{\alpha+1}) - f(u_\alpha) \over
    \delta}\, -f'(u)\Bigr|\\
  &=\Bigl|\sum_{j=1}^N {f(v_j) - f(v_{j-1}) \over \delta} \, -f'(u)\Bigr| \\
  &=\Bigl| \sum_{j=1}^N {f'(w_j) (v_j - v_{j-1}) \over \delta}\, -f'(u)\Bigr| \\
  &=\Bigl| \sum_{j=1}^N {v_j - v_{j-1}\over \delta}\,( f'(w_j) -f'(u))\Bigr| \\
  &\le \Lip(f') \sum_{j=1}^N {\left| v_j - v_{j-1}\right| \over
    \delta}\,\left| w_j -u \right| \\
  &\le \Lip(f') \delta\,,
\end{align*}
where each $w_j$, $j=1,\ldots, N$, is a suitable element in the
interval $(v_{j-1}, v_j)$. Passing to the essential supremum over $u$,
the proof is complete.$~~\diamondsuit$

\medskip

\begin{remark}\label{rem:homog_improved}
 For fixed $\nu\in\N$, the reconstructed flux $f_\nu$ could fail to 
 capture some inflection points of the real flux $f$, if e.g. the solution
 to the Riemann problem with datum~\eqref{eq:riemdata} consists
 of a single shock between $u_h$ and $u_{h+1}$ and the values of $f$ 
 in all inflection points in  $[u_h,u_{h+1}]$ are larger than 
 $\max\{f(u_h),f(u_{h+1})\}$.
 A possible way to overcome this intrinsic limitation could be to
 first fix $\nu\in\N$ and follow the proof of Theorem~\ref{thm:1},
 then to repeat the construction with $\nu'>\nu$ only in the subintervals
 of $[u_*,u^*]$ where the solutions corresponding to initial data 
 $u_o^h$ as in~\eqref{eq:riemdata}, with $h=0,\ldots, 2^\nu-1$,
 contained shocks. In this way, we focus our efforts on the intervals
 where the previous procedure might have been inaccurate, obtaining
 additional precision with a smaller number of tests, since not all the 
 states $u_k=u_*+2^{-\nu'} k$ get tested again.
\end{remark}

\subsection*{Proof of Theorem~\ref{thm:2}}
Fix a compact interval $J$ and denote by $\lambda=\max_{[u_1,u_2]}
|f'(u)|$, and set
\begin{equation*}
  I=\,(\min J-\lambda T,\max J+\lambda T)\,.
\end{equation*}
Now, consider the piecewise constant initial data given by
\begin{equation*}
  u_o(x)=\begin{cases}
    {\tilde u}, & \mbox{if } x\in I ,\\
    u_1, & \mbox{if } x\notin I,
  \end{cases}
\end{equation*}
for a fixed $\tilde u\in\,(u_1,u^m)$.
% we use the enlarged interval $I$ so that we have a natural bound to
% the wave location in $[x_M,\infty)$, at the end

By hyperbolicity there must exist a time $\tau>0$ small enough such
that the solution $u_{\rm obs}$ to the Cauchy problem
for~\eqref{eq:intrononh} with initial data $u(0,x)=u_o(x)$ can be
obtained, up to time $\tau$, by simply piecing together the solutions
to the Riemann problems with data
\begin{equation*}
  \begin{cases}
    u_1, & \mbox{if } x<\min I, \\
    {\tilde u}, & \mbox{if } x>\min I,
  \end{cases}\,\qquad\qquad \qquad
  \begin{cases}
    {\tilde u}, & \mbox{if } x <\max I, \\
    u_1, & \mbox{if } x>\max I.
  \end{cases}
\end{equation*}
Since it is not restrictive to assume $\tau\le T$, by the observabilty
assumption, both solutions to the Riemann problems above are
observable in $(0,\tau]$, and hence the whole solution to the Cauchy
problem with data $u_o$ is observable, up to time $\tau$.

Relying on the explicit construction of the solutions to Riemann
problems for~\eqref{eq:intrononh}, presented in~\cite{KlingRisebro}
and briefly sketched in the Appendix, we can also give a better a
priori description of the observed solution $u_{\rm obs}$ in
$[0,\tau]$. Indeed, the particular choice of initial data $u_o$ which
is constant in $I$, implies that any Lax wave present in $u_{\rm obs}$
for $x\in J$ must have been generated by a discontinuity in the flux
function $k$.

Moreover, the choice of a constant value $\tilde u<u^m$ in $u_o$
ensures that, at each discontinuity point $\xi$ for $k$, the solution
$u_{\rm obs}$ contains not only a stationary jump located at $x=\xi$,
but also a shock $u$-wave or a centered rarefaction $u$-wave with
positive speed (here and in the following $u$-waves are Lax waves with
constant values of $k$, see again the Appendix). Indeed, in a
neighborhood of $x=\xi$ the conservation law~\eqref{eq:intrononh} is
equivalent to a Riemann problem for the auxiliary
system~\eqref{eq:aux_sys} in the unknowns $(k,u)$ with initial data
\begin{equation*}
  \begin{cases}
    (\tilde u,k^\ell), & \mbox{ if }x<\xi,\\
    (\tilde u,k^r), & \mbox{ if }x>\xi,
  \end{cases}
\end{equation*}
for $k^\ell=k(\xi-)$ and $k^r=k(\xi+)$. Hence, the structure of the
solution can be deduced by the construction of the Riemann solver
for~\eqref{eq:aux_sys} with $u^\ell=u^r=\tilde u<u^m$ (see the
explicit description of the Riemann solver given in the Appendix, in
particular Cases~1 and~3).

Recalling that $k$ has a finite number of jumps in $J$, by {\bf (H1)},
there must be a time $\tau'>0$ such that the solution $u_{\rm obs}$ in
$(0,\tau']$ is obtained by piecing together the solutions of the
Riemann problems for the auxiliary system~\eqref{eq:aux_sys} at jumps
of the function $k$. In other words, at time $\tau'$ no interaction
between waves generated in $J$ has occurred yet. Without loss of
generality, we can assume that $\tau'\le \tau$.

We then observe the solution $u_{\rm obs}$ to~\eqref{eq:intrononh},
with initial data $u(0,x)=u_o(x)$, at times $t=\tau'/2$ and
$t=\tau'$. Denote by $S=\{x_1,\ldots, x_M\}$ the set (possibly empty)
of points $\xi\in J$ such that
\begin{equation*}
  u_{\rm obs}(\tau'/2, \xi-)=u_{\rm obs}(\tau', \xi-)\neq u_{\rm obs}(\tau',
  \xi+)=u_{\rm obs}(\tau'/2, \xi+)\,.
\end{equation*}
In other words, $S$ is the set of stationary jumps in the solution
$u_{\rm obs}$ and represents exactly the set of points of
discontinuity for the flux $k(x)$. Notice that $k(x)$ can have no
other jumps in $J$, because each jump in $k$ generates a stationary
discontinuity. This, in particular, implies that, if we find the
values of $k$ in the intervals $(x_{\alpha-1},x_\alpha)$, then we have
found exactly the correct function $k$ which produces $u_{\rm obs}$,
and~\eqref{eq:reconstr_2} is satisfied.

We introduce the notation
\begin{equation*}
  \kappa_o=k(x_1-)\,,\qquad\qquad\kappa_\alpha=k(x_\alpha+)\,,~~~~\alpha=1,\ldots,M\,,
\end{equation*}
for the values attained by $k(x)$ in $J$. From the admissibility of
the jumps in $u_{\rm obs}$, we deduce that for all
$\alpha\in\{1,\ldots, M\}$ one must have
\begin{equation*}
  \kappa_{\alpha-1} f(u_{\rm obs}(\tau', x_\alpha-))=\kappa_\alpha
  f(u_{\rm obs}(\tau', x_\alpha+))\,,
\end{equation*}
which is a set of $M$ equations in the $M+1$ unknowns
$\kappa_o,\ldots,\kappa_M$. To close the system we now need to find at
least one of the $\kappa_\alpha$. Indeed, if we can exactly identify
one of the unknowns, then the system above becomes a system in $M$
variables and $M$ unknowns, which can be solved because of the choice
$\tilde u>u_1$ in the initial data, which implies $f(u_{\rm
  obs}(\tau', x_\alpha\pm))>0$ for all $\alpha=1,\ldots,M$.

So we concentrate our attention on the interval $[x_M,\sup I]$ and we
define
\begin{equation*}
  y=\sup\big\{\xi>x_M~;~u_{\rm obs}(\tau',\xi)> \tilde u\big\}\,.
\end{equation*}
Having observed that in $[x_M,\sup I]$ there is a Lax $u$-wave
propagating with positive speed, $y$ is well defined and satisfies
$y\le\sup I$.
% any jump starting at $x_M\le\sup I$ can be located at most at
% $x_M+\lambda\tau'\le x_M+\lambda T=\sup I$ at time $\tau'$!
In particular, $y$ is the location, at time $t=\tau'$, of a $u$-wave
which got generated at $(t,x)=(0,x_M)$ from the discontinuity in $k$
and which is now moving away from $x_M$. The speed $\sigma$ of this
wave can be simply computed as
\begin{equation*}
  \sigma=\,{y-x_M\over \tau'}\,.
\end{equation*}
We want to find the value $\kappa_M$ from the speed of this $u$-wave
traveling in $[x_M,\sup I]$, so that the reconstruction of $k(x)$ is
complete and so is the proof.

There are two cases, depending on whether $u_{\rm obs}(\tau',\dott)$
is continuous or discontinuous at $y$. For ease of notation, define
\begin{equation*}
  u(y-):= u_{\rm obs}(\tau', y-)\,, \qquad\qquad u(y+):=
  u_{\rm obs}(\tau',y+)=\tilde u\,.
\end{equation*}
If $u(y-)=u(y+)$, then the $u$-wave is a centered rarefaction. In this
case,
\begin{equation*}
  \sigma=\kappa_M f'(\tilde u)\qquad \Longrightarrow \qquad
  \kappa_M=\,{\sigma\over f'(\tilde u)}\,,
\end{equation*}
which is well defined because $\tilde u<u^m$.

If $u(y-)\ne u(y+)$, then the $u$-wave is a shock and $u(y-)=u_{\rm
  obs}(\tau',x_M+)$. In this case,
\begin{equation*}
  \sigma=\kappa_M\, {f(\tilde u)-f(u(y-))\over \tilde u-u(y-)}\, \qquad
  \Longrightarrow \qquad \kappa_M=\sigma\,{\tilde u-u(y-)\over f(\tilde
    u)-f(u(y-))}\,,
\end{equation*}
which is well defined because $\sigma>0$ implies $f(\tilde u)\ne
f(u(y-))$.~~$\diamondsuit$

\subsection*{Proof of Theorem~\ref{thm:3}}
Most of the proof of Theorem~\ref{thm:3} follows from a series of
lemmas. The basic idea is that, since we assume only partial
observability for the solutions of every Cauchy problem, we
have to choose the initial data $u_o(\dott)$ in $(-\infty,a)$, so that
the observed solution gives enough data to reconstruct both the values
attained by $k(x)$ and the locations of its discontinuities inside
$(a,b)$.

We start by recalling that, as in Remark~\ref{rem:stationary_data},
the assumption that the initial data in $[a,\infty)$ is a stationary
solution to~\eqref{eq:intrononh} with $k$ and $f$ satisfying {\bf
  (H1)}--{\bf (H3)} implies that $u_o(x)$ has the
form~\eqref{eq:stationary_data} for some constant $\omega\in [u_1,u_2]$,
and that there holds the relation
\begin{equation}\label{eq:boundary_cond_2}
  k_o f(u_o(a))=k_1 f(\omega)=k_o f(u_o(b))\,,
\end{equation}
between the known quantities $k_o,u_o(a),u_o(b)$ and the unknown ones
$k_1,\omega$. From the analysis performed in~\cite{KlingRisebro}, we
also know that entropy admissibility of the stationary jumps located
at $x=\xi_1$ and $x=\xi_2$ (see the so-called ``smallest jump''
admissibility condition~\eqref{eq:smallest_jump}), implies that either
$f'(u_o(a)) f'(u_o(b))>0$, i.e., $u_m\notin\, (u_o(a),u_o(b))$, or
$\omega=u_m$ and hence $k_1$ is immediately determined
by~\eqref{eq:boundary_cond_2}. However, we have no way to determine
from the observations in $\R\setminus\,(a,b)$ which case is occurring
or which precise value $\omega$ is attained.

We now study the Cauchy problem with carefully selected initial data. Fix a
positive real value $\tilde x$ and assume $u_o(\dott)|_{[a,\infty)}$
to be a given stationary solution to~\eqref{eq:intrononh}. Let $\bar v_a$ be the state in $[u_1,u_2]$ characterized as the unique solution to 
\begin{equation}\label{eq:ua_companion}
f(u_o(a))=f(\bar v_a)\,,\qquad\qquad f'(u_o(a))f'(\bar v_a)\leq 0\,,
\end{equation}
so that, in particular, $\bar v_a\in(u_1,u^m)$ if $u_o(a)\in(u^m,u_2)$, and $\bar v_a\in(u^m, u_2)$ if $u_o(a)\in(u_1,u^m)$.
Define
\begin{equation*}
  \tilde y :=
  \begin{cases}
    \displaystyle\,{\max\{u_o(a),v_a\}-u_1 \over f(u_o(a))}\,
      f'(u_1)(b-a+\tilde x),  \qquad & \mbox{if } u_o(a)\in(u_1,u_2)\,,\\
    \,0, \qquad& \mbox{otherwise}\,,
  \end{cases}
\end{equation*}
and a
piecewise constant function $v_o$ as follows
\begin{equation}\label{eq:dato6}
  v_o(x)=
  \begin{cases}
    u^m, & \mbox{if } x< a-\tilde x-\tilde y \,,\\
    u_1,  & \mbox{if } a-\tilde x- \tilde y<x<a-\tilde x\,,\\
    u_o(a), \qquad &  \mbox{if } a- \tilde x<x<a\,,\\
    u_o(x), \qquad & \mbox{if } x>a\,.
  \end{cases}
\end{equation}
The following lemma helps to understand the choice of the initial data $v_o(\dott)$.
Namely, we show that after some time, the corresponding solution is identically equal 
to the state $u_1$ in $[a,b]$, and only afterwards the real reconstruction 
procedure begins, with larger values of the state variable crossing the unobserved region.
While this two--steps procedure is essential to remove
the possible presence of the state $u^m$ in the obstructed region $[\xi_1,\xi_2]$, 
which would  prevent the passage of any further wave through that region, 
it also implies that the procedure might require a large time of observation to be completed
 if, e.g., $u_o(a)$ is close to $u_1$. 
More comments on this aspect can be found in Remark~\ref{rem:faster_proc}.
\begin{lemma}\label{lem:3-3} Assume that the conservation law~\eqref{eq:intrononh}
  satisfies {\bf (H1)}--{\bf (H3)} and that $f(u)$ is a known
  function.  Let $u(t,x)$ denote the solution to the Cauchy problem
  for~\eqref{eq:intrononh} with initial data $v_o$ given by~\eqref{eq:dato6}, and 
  assume that $u_o(a)=v_o(a)\neq u^m$.

  Then, either $u_o(\dott)=u_2$ somewhere in $(a,b)$ or, by setting
  \begin{equation}\label{eq:tildetau}
  \tilde\tau:=\inf\left\{s>0~;~u(s,b)=u_1\right\}\,,
  \end{equation}
  with $\tilde\tau<+\infty$ and
  $u(\tilde\tau,\dott)|_{[a,b]}\equiv u_1$ and there exist times
  $T_1$, $T_2>\tilde\tau$ such that
  \begin{equation*}
    u_1<u(T_1,b)<u^m<u(T_2, a)<u_2\,.
  \end{equation*}
\end{lemma}

In the next lemma, we present a sufficient condition for finding a
unique solution to the inverse problem with prescribed stationary
initial data in $(a,\infty)$.

\begin{lemma}\label{lem:3-2} Assume that the conservation
  law~\eqref{eq:intrononh} satisfies {\bf (H1)}--{\bf (H3)}, that
  $f(u)$ is a known function and that the solution $u_{\rm obs}(t,x)$
  to the Cauchy problem for~\eqref{eq:intrononh} with initial data
  $v_o$ in~\eqref{eq:dato6} is partially observable in  $[0,T]\times\big(\R\setminus \,(a,b)\big)$.
  Then the following holds: if there exists $\tilde\tau\in[0,T]$ such that 
  $u_{\rm obs}(\tilde\tau,\dott)|_{[a,b]}\equiv u_1$ and if there exist 
  $T_1,T_2\in\,(\tilde\tau,T)$ such that
  \begin{equation}\label{eq:cond_stat}
  u_1=u_{\rm obs}(\tilde\tau,b)<u_{\rm obs}(T_1,b)<u^m<u_{\rm obs}(T_2, a)<u_2\,,
  \end{equation}
  then there exists a unique choice of $(k_1,\xi_1,\xi_2)$ such that,
  denoting with $u_{(k_1,\xi_1,\xi_2)}$ the solution to the Cauchy
  problem for~\eqref{eq:intrononh} with initial data $v_o$
  in~\eqref{eq:dato6} and with $k(x)$ given by~\eqref{eq:disck}, there
  holds
  $$
  u_{(k_1,\xi_1,\xi_2)}(t,x)=u_{\rm obs}(t,x), \quad
  (t,x)\in[0,T]\times(\R\setminus\,(a,b))\,.
  $$
\end{lemma}
The combination of the previous results simplifies the proof of
Theorem~\ref{thm:3}.

\medskip

\n {\it Proof of Theorem~\ref{thm:3}.} We claim that under the assumptions 
of the theorem, we have $u_o(a)\neq u^m$. Indeed, we are assuming 
that $u_o(\dott)$ is a stationary solution in $(a,b)$, and hence $u_o(a)=u^m$ 
would imply $u_o(x)\equiv u^m$ and $k(x)\equiv k_o$ on $[a,\infty)$. But this 
contradicts the assumption {\bf (H3)}, and hence it is not possible.

Since $u_o(a)\neq u^m$, we can choose to observe the solution 
corresponding to the initial data $v_o$ in~\eqref{eq:dato6} and 
combine Lemma~\ref{lem:3-3} and Lemma~\ref{lem:3-2} to conclude. 
Indeed, by Lemma~\ref{lem:3-3} we know that in finite time $t=\tilde\tau$, 
given by~\eqref{eq:tildetau}, the solution is constantly equal to $u_1$ in the 
unobservable region $(a,b)$. Moreover, there exist $T_1,T_2>\tilde\tau$ 
such that~\eqref{eq:cond_stat} holds, and hence Lemma~\ref{lem:3-2} ensures
the existence of a unique triple $(k_1,\xi_1,\xi_2)$ giving a solution
which satisfies~\eqref{eq:reconstr_3}. This concludes the
proof.~~$\diamondsuit$

\begin{remark} We observe that in the proof of Theorem~\ref{thm:3} 
we exclude the possibility of $u_o(a)=u^m$. It is clear that in that case the 
reconstruction is actually trivial: thanks to the assumption on $u_o(x)|_{[a,\infty)}$ 
being a stationary solution, the only possible flux function $k(x)$
has no jumps and it is constantly equal to $k_o$.  Such a reconstructed flux 
is excluded from the proof just because it does not satisfy the assumption {\bf (H3)}.
\end{remark}

\medskip

It remains now to prove Lemma~\ref{lem:3-3} and
Lemma~\ref{lem:3-2}.

\medskip

\n {\it Proof of Lemma~\ref{lem:3-3}.} 
In terms
of~\cite{KlingRisebro}, instead of~\eqref{eq:intrononh} we can study
the auxiliary system~\eqref{eq:aux_sys} for the unknowns $(u,k)$. In
this context, when dealing with piecewise constant initial data
like~\eqref{eq:dato6} we call $k$-wave
(resp.~$u$-wave) any Lax elementary wave, i.e., shock waves or
centered rarefaction waves, for the variable $k$ (resp.~$u$). It is
known (see Theorem~\ref{thm:KR}) that to each choice
$(k_1,\xi_1,\xi_2)$ in $(0,k_o]\times[a,b]\times[a,b]$, there
corresponds a unique entropy solution $u_{(k_1,\xi_1,\xi_2)}$ to the
Cauchy problem with initial data $v_o\in \BV(\R)$.

The choice of the initial data~\eqref{eq:dato6}, allows us to write explicitly 
the solution for small times (see the description of the Riemann solver
for~\eqref{eq:intrononh} in the Appendix). We focus our attention first to the 
case $u_o(a)>u_1$, so that $\tilde y>0$. Here, the solution
to~\eqref{eq:intrononh}, \eqref{eq:dato6} consists of a centered
rarefaction $u$-wave, starting at $a-\tilde x-\tilde y$ and evolving with characteristic 
speeds in $[0,k_o f'(u_1)]$, followed by a shock $u$-wave, starting at $a-\tilde x$ 
and traveling with speed
\begin{equation}\label{eq:shock_speed}
  \sigma=k_o\,{f(u_o(a)) - f(u_1)\over u_o(a)-u_1}\,=k_o\,{f(u_o(a))\over u_o(a)-u_1}\,,
\end{equation}
and by the stationary solution
$u_o(x)|_{[a,\infty)}$. And the structure of the solution is preserved at least  
as long as the shock $u$-wave remains in
$(-\infty,a)$. Notice that in the case under consideration the shock has strictly 
positive speed $\sigma$, because $u_o(a)>u_1$, and that we can write for all
$x\in\R$ and $t\in\left[0,{\tilde x \over \sigma}\,\right]$
\begin{equation}\label{eq:solution}
  u(t,x)=
  \begin{cases}
    u^m, & \mbox{if } x< a-\tilde x-\tilde y,\\
    \eta(x), & \mbox{if } a-\tilde x-\tilde y< x< a - \tilde x-\tilde y + \lambda_a t,\\
    u_1 & \mbox{if } a - \tilde x-\tilde y + \lambda_a t< x < a - \tilde x +\sigma t,\\
    u_o(a), & \mbox{if } a - \tilde x +\sigma t < x <a,\\
    u_o(x), & \mbox{if } x>a,
  \end{cases}
\end{equation}
where  $\lambda_a:=k_o f'(u_1)$ and $\eta(x)$ is the unique value such that 
$k_o f'(\eta(x))={x - (a  - \tilde x-\tilde y)\over t}$.

We prove that $\tilde\tau<+\infty$ and that $u(\tilde\tau,\dott)|_{[a,b]}\equiv u_1$, thanks 
to the choice of $\tilde y$. Indeed, the shock $u$-wave started at $a-\tilde x$ and 
traveling with speed $\sigma$ will eventually reduce its speed when 
it interacts with jumps of $k$, but it will always move with a speed 
$\sigma'\in[{k_o f(u_o(a)) \over M-u_1}\,,\,{k_o f(u_o(a)) \over m-u_1}]$, where
$$
m:= \min\big\{u_o(a),\bar v_a\big\}\,,\qquad M:=
\max\big\{u_o(a),\bar v_a\big\}\,,
$$ 
and $\bar v_a$ is the state characterized by~\eqref{eq:ua_companion}, as in the 
definition of $\tilde y$. This immediately implies that the wave will reach $x=b$ 
at most in time
$$
{M-u_1 \over k_o f(u_o(a))}\,(b-a+\tilde x)<+\infty
$$
and that such a time gives an upper bound to $\tilde\tau$.
Moreover, the choice of $\tilde y$ now implies that the rarefaction $u$-wave
generated by the jump at $x=a-\tilde x-\tilde y$ is still traveling in
$[a-\tilde x-\tilde y,a-\tilde x]$ when the shock emerges at $x=b$.
This implies that $u_{\rm obs}(\tilde \tau,\dott)$ is a stationary
solution for~\eqref{eq:intrononh} in $(a,b)$ with
$$
u_{\rm obs}(\tilde \tau,a)=u_{\rm obs}(\tilde \tau,b)=u_1<u_m\,.
$$
Observing that $u_{\rm obs}(\tilde\tau,\dott)|_{[b,\infty)}$ contains a
single shock wave traveling with positive speed, and therefore moving
away from the unobservable region $(a,b)$, this shock will not contribute anymore
to the values attained by the solution in $(a,b)$ for times $t\geq \tilde\tau$.

Consider now the case in which $u_o(a)=u_1$, and hence $\tilde y=0$ in~\eqref{eq:dato6}.
In this case, the shock $u$-wave is not present at all and the solution $u(t,x)$ 
in~\eqref{eq:solution} attains value $u_o(a)=u_1$ for $x\in[a - \tilde x-\tilde y + \lambda_a t,a]$ 
and times  $t\in\left[0,{\tilde x \over \lambda_a}\,\right]$, noticing that $\lambda_a>0$ thanks
to the assumption $u_o(a)\neq u^m$. 
Since the assumptions $u_o(\dott)<u_2$ and $u_o$ stationary solution in $[a,b]$ imply 
that $u_o\equiv u_1$, then we can conclude that $\tilde\tau = 0$ and that 
$u(\tilde\tau,\dott)|_{[a,b]}\equiv u_1$ as before.

This completes the first step of the procedure, needed to remove the
possible presence of congested regions. In the rest of the proof, we analyze
the evolution of the solution for times larger than $\tilde \tau$ in
order to reconstruct $k$ in $(a,b)$.

\smallskip

For times $t\geq\tilde\tau$, the rarefaction $u$-wave approaches the obstructed region 
and eventually reaches $x=\xi_1$ at time 
$t_{\xi_1}={\xi_1-a+\tilde x+\tilde y\over \lambda_a}$, which is unknown since $\xi_1\in[a,b]$
is unknown. Since $u(t_{\xi_1},\xi_1+)=u_1$ from the previous analysis,
after the interaction between the rarefaction wave and the stationary jump in $k$ at $x_1$, 
part of the wave simply passes through the obstruction. 
The result for $x>\xi_1$ would
then be a new rarefaction $u$-wave, propagating with a smaller
characteristic speed.
% and connecting the states $u_1$ and $u^m$, where
%$w$ is the only state in $[u_1,u^m)$ such that $k_o f(u_o(a))= k_1
%f(w)$. 
This new centered rarefaction $u$-wave is going to pass through
$x=\xi_2$ at some later time $t_{\xi_2}$, and it keeps propagating towards $x=b$, 
because we also have that $u(t_{\xi_2},\xi_2+)=u_1$.
%thanks to the assumption $u_o(b)<u^m$ in the
%definition~\eqref{eq:dato4} of $v_o$. 

Notice that, for times $t\geq t_{\xi_1}$, the value $u(t,\xi_1-)$
increases due to the incoming rarefaction wave and the solution
$u(t,x)$ for $x\in[\xi_1,\xi_2]$ will be a smooth profile
corresponding to a rarefaction $u$-wave joining $u_1$ with the value
$u(t,\xi_1+)$ characterized by being the only state in $(u_1,u^m)$
with the property
$$
k_1 f(u(t,\xi_1+))=k_o f(u(t,\xi_1-))\,.
$$
Since $k_1<k_o$, the region $[\xi_1,\xi_2]$ becomes congested before
the whole original $u$-rarefaction can pass through $x=\xi_1$. More
precisely, setting $w:=u(t,\xi_1-)$ the state for which
$u(t,\xi_1+)=u^m$, then $w$ is the maximal value of the
conserved quantity that the obstructed region $[\xi_1,\xi_2]$ can
accept. However, due to the continuous arrival of larger states from
the left side, a shock $u$-wave appears at $x=\xi_1-$ and travels back
towards $x=a$ with negative speed. Notice that along such a
``reflected'' discontinuity, the right state is always given by
$w'\in(u^m,u_2)$ such that $k_o f(w')= k_1 f(u^m)=k_o f(w)$.

We sum up the discussion so far: Due to the propagation of the smaller
states of the rarefaction wave, we find $T_1>\tilde\tau$ such that $u(T_1,
b)>u(\tilde\tau, b)=u_1$; due to the reflected shock which emerges at $x=\xi_1$,
there exists $T_2>\tilde\tau$ such that $u(T_2,a)>u^m$. Therefore, the lemma is
proved.~~$\diamondsuit$

\medskip

\n {\it Proof of Lemma~\ref{lem:3-2}.}  Set
\begin{equation*}
  \tau_a:= \inf \{t > \tilde\tau~;~ u_{\rm obs}(t, a) > u^m\}\,,
\end{equation*}
and
\begin{equation*}
  \tau_b := \inf \{ t > \tilde\tau~;~ u_{\rm obs}(t, b) > u_1\}\,.
\end{equation*}
These are known values, thanks to the partial observability assumption and we have $\tilde\tau<\tau_a\le T_2$ and
$\tilde\tau<\tau_b<T_1$. 

\n Indeed, the description of the Riemann solver
for~\eqref{eq:intrononh} given in the Appendix, implies that in the 
case $u_o(a)>u_1$, so that $\tilde y>0$, for small
positive times $u_{\rm obs}$ consists of a centered
rarefaction $u$-wave, starting at $a-\tilde x-\tilde y$ and evolving 
with characteristic 
speeds in $[0,k_o f'(u_1)]$, followed by a shock $u$-wave, starting 
at $a-\tilde x$ 
and traveling with speed $\sigma$, as in~\eqref{eq:shock_speed}, 
and by the stationary solution
$u_o(x)|_{[a,\infty)}$. And in the case  $u_o(a)=u_1$ the structure is similar but 
without the shock $u$-wave.

\n Since we are assuming that  $u(\tilde\tau,\dott)|_{[a,b]}\equiv u_1$, 
this means that at time $\tilde \tau$ the shock wave has already passed 
through the whole unobservable region and the rarefaction wave has not reached it
yet. Then, it follows that $\tau_b>\tilde\tau$ is the first time when the centered 
rarefaction appears at the end of the unobservable region, while $\tau_a>\tilde\tau$ 
is the first time when a shock is reflected by the discontinuities of $k$ inside $(a,b)$
back towards $x=a$. Moreover, since $k$ has the form~\eqref{eq:disck},
we also know that this shock $u$-wave emerging at time $\tau_a$
originated at $x=\xi_1$, when the rarefaction $u$-wave above interacted
with the stationary jump of $k(x)$ and the state at $x=\xi_1+$ reached
the value $u^m$. For later use, let us define $\tau_o\geq \tilde\tau$ the first time when
the rarefaction $u$-wave originated at $a-\tilde x-\tilde y$ reaches $x=a$. Thanks to the 
partial observability assumption, $\tau_o$ can be considered a known value.

Let now $v := u(\tau_a, a-)$ and $w' := u(\tau_a, a+)$ be the states
separated by the shock wave emerged at $x=a$. By the Rankine--Hugoniot
conditions at the generating point of the shock $u$-wave, there must
hold
\begin{equation*}
  k_o f(w')=k_1 f(u^m)\,,
\end{equation*}
which in turn implies
\begin{equation}\label{eq:k1}
  k_1 = \,{k_o f(w') \over f(u^m)}\,.
\end{equation}
Focusing our attention on the evolution of the rarefaction $u$-wave in
$(a,b)$, we know that $\tau_b-\tau_o$ must be equal to
\begin{equation*} {\xi_1 - a \over k_o f'(u_1)} + {\xi_2 - \xi_1
    \over k_1f'(\omega)} + {b - \xi_2 \over k_o f'(u_1)}\,,
\end{equation*}
and $\omega$ is the unique solution in $(u_1,u_m]$ of
\begin{equation}\label{eq:omega}
  k_o f(u_1)=k_1 f(\omega).
\end{equation}
Indeed, the wave must have traveled with speed $k_o f'(u_1)$ in
$[a,\xi_1)$, with speed $k_1f'(\omega)$ in $(\xi_1,\xi_2)$, and
again with speed $k_o f'(u_1)$ in $(\xi_2,b]$. Note that $\omega$
is now known from~\eqref{eq:omega}, because $f$ is a known function
and $k_1$ has been already found in~\eqref{eq:k1}. Therefore, we
obtain
\begin{equation}\label{eq:ampiezza}
  \xi_2 - \xi_1 =\,{k_o k_1 f'(\omega) \over k_o f'(u_1) - k_1 f'(\omega)} \left[ {b - a\over k_o}\, -(\tau_b-\tau_o) f'(u_1) \right]\,,
\end{equation}
where all quantities appearing at the right-hand side are known.

We want to use $\tau_a$ and the states $v$, $w'$ observed in
$(\tau_a,a)$ to determine $\xi_1$. Let $w\in\,[u_1,u^m)$ be the unique
solution of $f(w)=f(w')$. We know from the structure of the Riemann
solver that the shock $u$-wave separating $v$ and $w'$ originated at
$x=\xi_1$ when the rarefaction wave traveling with speed $k_o f'(w)$
interacted with the stationary jump of $k(x)$ and $u_{\rm obs}$ at
$x=\xi_1+$ reached $u^m$. Then, we can conclude that the interaction
at $x=\xi_1$ which generated the reflected shock occurred at the time
\begin{equation*}
  \bar\tau (\xi_1):= {\xi_1 - a + \tilde x + \tilde y \over k_o f'(w)}\leq \tau_a \,.
\end{equation*}

Notice that, due to the structure of the rarefaction wave, we have
that $u_1<w\leq v= u(\tau_a, a-)$.  Indeed, it is not possible to have
$w=u_1$ as the reflected state because in such a case we would have
$w'=u_2$ and
$$
k_1=\,{k_o f(u_2)\over f(u^m)}\,=0\,,
$$
which is not possible for a function $k$ satisfying {\bf (H3)}. On the other hand, the limit case 
$w=v$ happens when the reflection occurs at time $\tau_a$ and hence
it is equivalent to having $\xi_1=a$. Since $\xi_2$ is uniquely determined as well, by 
using~\eqref{eq:ampiezza}, the proof is complete.

It remains to consider the case $u_1<w<v$. This means we can assume $\xi_1>a$ and, 
hence, there hold both 
$\bar\tau(\xi_1)<\tau_a$ and $\bar\tau(\xi_1)>{\tilde x+\tilde y\over k_o f'(w)}$, because the latter is
the time at which the rarefaction front passes at $x=a$, before getting reflected.
Moreover, the wave observed in $x=a$ at $t=\tau_a$ is exactly the
(forward) generalized characteristic $\xi(t)$ associated to $u=u_{\rm
  obs}$, emanating from the point $(\bar\tau(\xi_1), \xi_1)$
(see~\cite{DafermosBook}). Due to the particular structure of our
problem, this curve can be found as the solution of the backward
Cauchy problem for
\begin{equation}\label{eq:genchar}
  \dot\xi(t)=k_o \,{f\left(u(t,\xi(t)-)\right) - f\left(w'\right)
    \over u(t,\xi(t)-) - w'}\,, 
\end{equation}
with data
\begin{equation}\label{eq:data_genchar}
  \xi(\tau_a)=a\,.
\end{equation}
Thanks to the regularity of $f$ in $[u_1,u_2]$ and of $u$ in
$\Omega=(-\infty,\xi_1)$, the
problem~\eqref{eq:genchar}--\eqref{eq:data_genchar} has a unique
Carath\'eodory solution defined in $(\bar\tau(\xi_1),\tau_a]$, since
we know that only at $t=\bar\tau(\xi_1)$ the solution reaches the boundary
$\partial\Omega$.  Therefore, $\xi_1$ satisfies the relation:
\begin{equation}\label{eq:xi_implicit}
  \xi_1=a - \int_{\bar\tau(\xi_1)}^{\tau_a} \dot \xi(t)\,dt=
  a - k_o \int_{\bar\tau(\xi_1)}^{\tau_a}{f\left(u(t,\xi(t)-)\right)
    - f\left(w'\right)  \over u(t,\xi(t)-) - w'}\,dt\,. 
\end{equation}
Setting
\begin{equation*}
  \chi(\xi_1):= \xi_1 - a + k_o
  \int_{\bar\tau(\xi_1)}^{\tau_a}{f\left(u(t,\xi(t)-)\right) -
    f\left(w'\right) \over u(t,\xi(t)-) - w'}\,dt\,,
\end{equation*}
we can combine
\begin{align*}
  {d\chi\over d\xi_1} &= 1 - k_o \,{f\left(u(\bar\tau,\xi_1)\right) -
    f\left(w'\right) \over u(\bar\tau,\xi_1) -
    w'} {d\bar\tau\over d\xi_1} \\
  &= 1 - {1\over f'(w)} {f\left(w\right) - f\left(w'\right) \over w -
    w'} = 1>0\,,
\end{align*}
with $\chi(a)< 0$ and $\chi(a-\tilde x-\tilde y+k_o f'(w) \tau_a)=
k_o f'(w)\left(\tau_a-{\tilde x+\tilde y\over k_o f'(w)}\right)>0$, to
conclude that there exists a unique value $\xi_1$ such that
$\chi(\xi_1)=0$, i.e., a unique location $\xi_1$ where the reflected
shock has been generated. Finally, using~\eqref{eq:ampiezza},
$\xi_2$ is uniquely determined as well and the proof is
complete.~~$\diamondsuit$

\medskip

\begin{remark} It is worth noticing that given $f$, the expression for
  $\xi_1$ can be explicitly obtained from~\eqref{eq:xi_implicit}. To
  fix ideas, let $[u_1,u_2]=[0,1]$ and $f(u)=u(1-u)$, as in
  Example~\ref{ex:highway}. Then the ordinary differential equation
  solved by $\xi_1$ reduces to
  \begin{equation*}
    \dot \xi(t)= k_o\left(1 - w' - u(t,\xi(t)-)\right)
  \end{equation*}
  which implies that~\eqref{eq:xi_implicit} can be written in the form
  \begin{equation*}
    \xi_1=a - k_o(\tau_a-\bar\tau) (1-w') - k_o \int_{\bar\tau}^{\tau_a}
    u(t,\xi(t)-) \,dt.
  \end{equation*}
  Since in this case, the solution $u_{\rm obs}$ is given
  by~\eqref{eq:solution} with $\eta(x)={1\over 2} - {\xi_1 - a +
    \tilde x \over 2k_o t}$, the integral can be computed explicitly
  and $\xi_1$ can be retrieved as a root of a polynomial of degree
  three.
\end{remark}

\begin{remark}\label{rem:faster_proc}
 The choice of the initial data  $v_o(\dott)$ in~\eqref{eq:dato6} for the proof of Theorem~\ref{thm:3} 
 needs a few comments. With such a choice, the reconstruction procedure consists in 
 ``emptying'' the unobservable region $(a,b)$ before starting to send new waves that allow
 to  identify exactly the location and size of the obstruction.
 The first part of the procedure cannot be avoided when $\max\{u_o(a), u_o(b)\}> u^m$,
 because in this case no rarefaction wave can pass through
 the congested part of $(a,b)$ to collect the information needed for the reconstruction. 
 However, this makes the process slower whenever no congested region is present.

 An alternative choice when $\max\{u_o(a), u_o(b)\}< u^m$ is the following.
 The assumption that $u_o(\dott)|_{[a,\infty)}$ is a stationary solution to~\eqref{eq:intrononh}, implies that
 the flux function $k(\dott)$ attains a value $k_1\in[k_o\,{f(u_o(a))\over f(u^m)},k_o)$ in $[\xi_1,\xi_2]$.
 If we knew that $k_1>k_o\,{f(u_o(a))\over f(u^m)}$, then a more effective choice of the initial 
 data would be
 \begin{equation*}
  w_o(x)=
  \begin{cases}
    u^m,  \qquad & \mbox{if } x<a-\tilde x\,,\\
    u_o(a), \qquad& \mbox{if } a-\tilde x<x<a,\\
    u_o(x), \qquad& \mbox{if } x>a\,,
  \end{cases}
  \end{equation*}
 for any choice of $\tilde x>0$. In the solution to the Cauchy problem for~\eqref{eq:intrononh} 
 with initial data $w_o(\dott)$ there is no shock wave emptying the 
 unobservable region, but only the rarefaction wave connecting 
 the states $u_o(a)$ and $u^m$. Hence, it would still be possible
 to proceed as in the proof of Lemmas~\ref{lem:3-3}--\ref{lem:3-2} 
 and to find a unique flux function $k(\dott)$
 with the properties required in Theorem~\ref{thm:3}. In addition, the process could be completed 
 in a shorter time.
 
 \n The problem is that a priori we cannot exclude that  $k_1=k_o\,{f(u_o(a))\over f(u^m)}$, 
 or equivalently that $u_o(x)=u^m$ for $x\in[\xi_1,\xi_2]$, and in this case
 no wave in the solution would pass through the congested region $[\xi_1,\xi_2]$.
 Thus, repeating the previous reconstruction procedure would only give $k_1$ and $\xi_1$, 
 but not $\xi_2$.
 
 A way to combine the best aspects of both approaches is to use $w_o(\dott)$ as 
 initial data and wait to see if at some time $\tau>0$ a shock $u$-wave appears at $x=a$, 
 separating the states
 $u_{\rm obs}(\tau, a-)$ and $u_{\rm obs}(\tau, a+)>u^m$ with
 $$
 f(u_{\rm obs}(\tau, a+))=f(u_o(a))\,.
 $$
 If this happens, then we realize a posteriori that the region $[\xi_1,\xi_2]$ was 
 originally congested. Therefore, relying on the fact that $u_{\rm obs}(\tau, \dott)|_{[a,\infty)}$ 
 is a stationary solution for~\eqref{eq:intrononh} with $u_{\rm obs}(\tau, a)>u^m$, 
 we can restart the procedure for times $t\geq\tau$ with the initial data 
 $v_o(\dott)$ given by~\eqref{eq:dato6} and complete the reconstruction process.
 
 \n On the other hand, if $f(u_{\rm obs}(\tau, a+))>f(u_o(a))$ or if there exists $\tau'>0$ 
 such that $u_{\rm obs}(\tau', b)>u_o(b)$, then we can deduce that no congested area was 
 present in $(a,b)$ at time $t=0$ and the initial data $w_o(\dott)$ will be sufficient to
 complete the reconstruction procedure.
\end{remark}

\subsection*{Proof of Theorem~\ref{thm:4}}
The proof of Theorem~\ref{thm:4} follows from two lemmas. First of
all, note that if $k_1\ne k_o$, then $u\equiv \bar u_o$ is not a
stationary solution for~\eqref{eq:intrononh}. Hence, at time $t=0+$
the jumps in $k(x)$ produce waves with non-zero speed at one or both
sides of each jump (see the Appendix). Our first lemma deals with the
evolution of the solution to~\eqref{eq:intrononh} when $u(0,x)\equiv
\bar u_o$ and $k$ is of the form~\eqref{eq:disck}.

\begin{lemma}\label{lem:4-1} Assume that the conservation
  law~\eqref{eq:intrononh} satisfies {\bf (H1)}--{\bf (H3)} and that
  $f(u)$ is a known function. Fix a constant $\bar u_o\in[u_1,u^m)$
  and denote by $u(t,x)$ the solution to the Cauchy problem
  for~\eqref{eq:intrononh} with $u(0,x)\equiv \bar u_o$ for
  $x\in\R$. Then, the following facts hold:
  \begin{description}
  \item[(i)] If $k_o f(\bar u_o)> k_1 f(u^m)$, then there exist
    $T_1,T_2>0$ such that
    \begin{equation}\label{eq:cond_nstat_1}
      u(T_1,b+)=\bar u_o> u(T_1, b-)\mbox{ ~~~and~~~ } u(T_2,a-)=\bar
      u_o< u(T_2,a+)\,. 
    \end{equation}
  \item[(ii)] If $k_o f(\bar u_o)\le k_1 f(u^m)$, then either
    $u(t,\dott)\equiv \bar u_o$ in $\R\setminus(a,b)$ for all $t>0$ or
    there exist $0<T_1<T_2$ such that
    \begin{equation}\label{eq:cond_nstat_2}
      u(T_1,b+)=\bar u_o> u(T_1, b-)\mbox{ ~~~and~~~ } u(T_2,b)> u(T_1,b-)\,.
    \end{equation}
    Moreover, the former case does not happen if $k_o f(\bar u_o)= k_1
    f(u^m)$.
  \end{description}
\end{lemma}

The second lemma gives sufficient conditions for the existence of a
triple $(k_1,\xi_1,\xi_2)$ such that the solution
to~\eqref{eq:intrononh} with constant initial data $u(0,x)\equiv \bar
u_o$ and flux $k(x)f(u)$, $k(x)$ being given by~\eqref{eq:disck},
coincides with the  observed solution $u_{\rm obs}$ in the
observable region.

\begin{lemma}\label{lem:4-2} Assume that the conservation
  law~\eqref{eq:intrononh} satisfies {\bf (H1)}--{\bf (H3)}, that
  $f(u)$ is a known function and that the solution $u_{\rm obs}(t,x)$
  to the Cauchy problem for~\eqref{eq:intrononh} with constant initial
  data $u(0,x)\equiv \bar u_o\in[u_1,u^m)$ is partially observable in
  $[0,T]\times\big(\R\setminus\,(a,b)\big)$.

  \begin{description}
  \item[(i)] If there exist $T_1,T_2\in\,(0,T)$ such that
    \begin{equation}\label{eq:cond_nstat_1_thm}
      u_{\rm obs}(T_1,b+)=\bar u_o> u_{\rm obs}(T_1, b-)\mbox{ ~~~and~~~ }
      u_{\rm obs}(T_2,a-)=\bar u_o< u_{\rm obs}(T_2,a+)\,,  
    \end{equation}
    then there exists a unique choice of $(k_1,\xi_1,\xi_2)$ such that
    if $u_{(k_1,\xi_1,\xi_2)}$ denotes the solution of the Cauchy
    problem for~\eqref{eq:intrononh} with initial data $u(0,x)\equiv
    \bar u_o$ and $k(x)$ given by~\eqref{eq:disck}, we have that
    \begin{equation*}
      u_{(k_1,\xi_1,\xi_2)}(t,x)=u_{\rm obs}(t,x), \quad  (t,x)\in[0,T]\times(\R\setminus\,(a,b))\,.
    \end{equation*}
  \item[(ii)] If there exist $0<T_1<T_2<T$ such that
    \begin{equation}\label{eq:cond_nstat_2_thm}
      u_{\rm obs}(T_1,b+)=\bar u_o> u_{\rm obs}(T_1, b-)\mbox{ ~~~and~~~ }
      u_{\rm obs}(T_2,b)> u_{\rm obs}(T_1,b-)\,, 
    \end{equation}
    then there exists a choice of $(k_1,\xi_1,\xi_2)$ such that if
    $u_{(k_1,\xi_1,\xi_2)}$ denotes the solution of the Cauchy problem
    for~\eqref{eq:intrononh} with initial data $u(0,x)\equiv \bar u_o$
    and $k(x)$ given by~\eqref{eq:disck}, we have that
    \begin{equation*}
      u_{(k_1,\xi_1,\xi_2)}(t,x)=u_{\rm obs}(t,x), \quad
      (t,x)\in[0,T]\times(\R\setminus\,(a,b))\,.
    \end{equation*}
    Moreover, if
    \begin{equation*}
      \inf\left\{s\in~(T_1,T_2)~;~u(s,b)>u(T_1,b-)\right\}> T_1\,,
    \end{equation*}
    then the choice is also unique.
  \end{description}
\end{lemma}
Now the proof of Theorem~\ref{thm:4} is immediate.

\medskip

\n{\it Proof of Theorem~\ref{thm:4}.} Let $u_{\rm obs}(t,x)$ denote
the solution of the Cauchy problem for~\eqref{eq:intrononh} with
constant initial data $u(0,x)\equiv \bar u_o\in[u_1,u^m)$ and a flux
function $kf$ with $k$ given by~\eqref{eq:disck}. Even if we do not
know the values of $k_1,\xi_1,\xi_2$, we know that either $k_o f(\bar
u_o)> k_1 f(u^m)$ or $k_o f(\bar u_o)\le k_1 f(u^m)$. In the former
case, Lemma~\ref{lem:4-1} ensures that there exist $T_1,T_2>0$ such
that~\eqref{eq:cond_nstat_1} holds for $u_{\rm obs}$. Hence, we can
apply part {\it (i)} of Lemma~\ref{lem:4-2} to find the triple
$(k_1,\xi_1,\xi_2)$ which gives a solution
satisfying~\eqref{eq:reconstr_4}. Similarly, in the latter case,
Lemma~\ref{lem:4-1} ensures that either $u(t,\dott)\equiv \bar u_o$ in
$\R\setminus(a,b)$ for all $t>0$, or there exist $0<T_1<T_2$ such
that~\eqref{eq:cond_nstat_1} holds for $u_{\rm obs}$. In particular,
if $u(t,\dott)\not\equiv \bar u_o$ for all times $t>0$, part {\it
  (ii)} of Lemma~\ref{lem:4-2} gives $(k_1,\xi_1,\xi_2)$ such
that~\eqref{eq:reconstr_4} holds.

The uniqueness part follows from Lemma~\ref{lem:4-2} as well, under
the hypotheses of Theorem~\ref{thm:4}, completing the
proof.~~$\diamondsuit$

\medskip

It remains to prove Lemma~\ref{lem:4-1} and Lemma~\ref{lem:4-2}.

\medskip

\n {\it Proof of Lemma~\ref{lem:4-1}.} In terms
of~\cite{KlingRisebro}, we can study the Cauchy problem
for~\eqref{eq:intrononh} by studying the auxiliary
system~\eqref{eq:aux_sys} for the unknowns $(k,u)$. In this case, the
initial data $u(0,\dott)\equiv \bar u_o$ 
is written
\begin{equation}\label{eq:aux_initial}
  (u_{in},k_{in})=(u,k)(0,x) = 
  \begin{cases}
    (\bar u_o, k_o) & \mbox{if } x<\xi_1\,,\\
    (\bar u_o, k_1) & \mbox{if } \xi_1<x<\xi_2\,,\\
    (\bar u_o, k_o) & \mbox{if } x>\xi_2\,,\\
  \end{cases}
\end{equation}
for the resonant system~\eqref{eq:aux_sys}. In the following, we call
$u$-waves any Lax elementary wave for~\eqref{eq:aux_sys}, i.e., shock
waves or centered rarefaction waves, propagating with constant $k$,
and $k$-waves the stationary jumps between states
satisfying~\eqref{eq:rh_stationary} (see the Appendix). When
considering the Cauchy
problem~\eqref{eq:aux_sys}--\eqref{eq:aux_initial}, both jumps in the
initial data $(u_{in},k_{in})$ create for $t>0$ a stationary $k$-wave
and some $u$-waves. Namely, since $\xi_1< \xi_2$, for small times the
solution to the Cauchy problem is given by the juxtaposition of the
solutions to the Riemann problems in $x=\xi_1$ and $x=\xi_2$ and it
can be described, in terms of the Riemann solver described in the
Appendix, as follows:

\smallskip

\n {\it (i)} Assume $k_o f(\bar u_o)> k_1 f(u^m)$ and consider first
the jump in $x=\xi_1$. In this case, the Riemann problem is solved by
a shock $u$-wave with negative speed $\sigma^-$, a
stationary $k$-wave and a centered rarefaction $u$-wave  with
positive speeds. In terms of the original system, this means that at
$x=\xi_1+$ the variable passes from $\bar u_o$ to the larger $u^m$,
because there is more incoming $u$ than the obstructed region can
carry. Moreover, a shock wave appears in $x=\xi_1-$ and propagates
back towards $x=a$ with speed $\sigma^-$. In terms of the traffic flow
model presented in Example~\ref{ex:highway}, such a solution can be
interpreted as a queue of cars forming at $x=\xi_1-$ and traveling
back towards $x=a$, followed by a region of congested traffic in
$x=\xi_1+$ due to the continuous arrival of more cars than the
obstructed highway can carry.

Consider now the jump from $k_1$ to $k_o$ at $x=\xi_2$. The Riemann
problem for~\eqref{eq:aux_sys} is solved by a stationary $k$-wave
followed by a shock $u$-wave traveling with positive speed $\sigma^+$
towards $x=b$. In terms of the original system, this means that at
$x=\xi_2+$ a smaller value $u'\in [u_1,\bar u_o]$, such that
$k_of(u')=k_1 f(\bar u_o)$, emerges from the discontinuity and
therefore a shock $u$-wave between $u'$ and $\bar u_o$ forms. From the
point of view of Example~\ref{ex:highway}, this means that cars at
$x=\xi_2-$ have to slow down due to the obstruction and a region with
smaller car density $u'$ appears at $x=\xi_2+$.

As $t$ increases, the $u$-shock traveling with speed $\sigma^-$ simply
propagates in $(a,\xi_1)$ and reaches $x=a$ at time $T_2={a-\xi_1
  \over \sigma^-}>0$, as requested by the second part
of~\eqref{eq:cond_nstat_1}.

On the other hand, the rarefaction $u$-wave created at $x=\xi_1$
eventually interacts with the stationary $k$-wave in $x=\xi_2$ and
keeps propagating in $(\xi_2,b)$ as a rarefaction $u$-wave, but with
larger positive speeds.  In particular, this rarefaction $u$-wave in
$(\xi_2,b)$ is now separating the states in the interval $[u',\bar
u_o)$ and its front travels with speed $f'(u')$ larger than the
speed $\sigma^+$ of the shock $u$-wave generated at $t=0$ in
$x=\xi_2$, due to the entropy admissibility of the shock.  Hence, the
$u$-rarefaction could start interacting with the $u$-shock, before
they reach $x=b$. However, the $u$-shock cannot be completely canceled
by the $u$-rarefaction and, therefore, there must exist $T_1>0$ such
that the first part of~\eqref{eq:cond_nstat_1} is verified as well.

\smallskip

\n {\it (ii)} Assume now $k_o f(\bar u_o)\le k_1 f(u^m)$, and consider
first the jump at $x=\xi_1$. In this case, the Riemann problem is
solved by a stationary $k$-wave and a centered rarefaction $u$-wave
traveling with positive speed. In terms of the original system, this
means that the incoming quantity $\bar u_o$ from $x=\xi_1-$ does not
completely fill the region in $(\xi_1,\xi_2)$. As a consequence, $u$
only increases from $\bar u_o$ to a larger value $u''\le u^m$ such
that $k_of(\bar u_o)=k_1 f(u'')$.  In the terminology of
Example~\ref{ex:highway}, this means that the incoming cars do not
completely fill the road in $(\xi_1,\xi_2)$ and therefore no queue
appears at $x=\xi_1-$.

Considering the jump from $k_1$ to $k_o$ at $x=\xi_2$, the Riemann
problem for~\eqref{eq:aux_sys} is solved again by a stationary
$k$-wave followed by a shock $u$-wave traveling with positive speed
$\bar\sigma$ towards $x=b$. As before, in terms of the original
system, this means that at $x=\xi_2+$ a smaller value $u'\in [u_1,\bar
u_o]$, such that $k_of(u')=k_1 f(\bar u_o)$, emerges from the
discontinuity and therefore a shock $u$-wave between $u'$ and $\bar
u_o$ forms.

As $t$ increases, the rarefaction $u$-wave exiting $x=\xi_1$ will
eventually interact with the stationary $k$-wave in $x=\xi_2$ and will
keep propagating as a rarefaction $u$-wave, but with larger positive
speed in $(\xi_2,b)$. As before, this rarefaction $u$-wave is faster
than the shock $u$-wave traveling with speed $\bar\sigma$ and, hence,
the $u$-waves could start interacting before reaching $x=b$. This
interaction opens up two different scenarios:
\begin{itemize}
\item Either the whole interaction between rarefaction and shock takes
  place in $(\xi_2, b)$, resulting in a complete cancellation of the
  two waves (this is the case whenever
  \begin{equation*}
    {\xi_2-\xi_1 \over k_1 f'(\bar u_o)}\,+\,{b-\xi_2 \over k_o
      f'(u')}\,\le\,{b-\xi_2 \over \bar\sigma}\,,
  \end{equation*}
  where $\bar\sigma$ denotes the speed of the shock separating $u'$
  and $\bar u_o$).
\item Or there exists $\tau_1>0$ such that $u(\tau_1,b)<\bar u_o$.
\end{itemize}
In the former case, $u(t,\dott)|_{\R\setminus(a,b)}\equiv \bar u_o$
for all $t>0$. In the latter case, we simply set $T_1=\tau_1$ and
\begin{equation*}
  T_2=\inf\left\{s\geq T_1~;~ u(s,b-)>u(T_1,b)\right\} + \ve\,,
\end{equation*}
for any fixed $\ve>0$
small. %the \ve is needed in the case of an ongoing interaction when the waves emerge at x=b!
Observe that if $k_o f(\bar u_o)= k_1 f(u^m)$ the last case holds, and
hence $u''=u^m$, which completes the proof.~~$\diamondsuit$

\medskip

\n {\it Proof of Lemma~\ref{lem:4-2}.} By assuming that $f$ and $k$
satisfy {\bf (H1)}--{\bf (H3)}, one immediately obtains some
properties of the solution $u_{\rm obs}$ to~\eqref{eq:intrononh} with
constant initial data $u(0,x)\equiv \bar u_o$. In particular, from the
description of the Riemann solver for~\eqref{eq:intrononh} given in
the Appendix, one can see that $k_1< k_o$ in~\eqref{eq:disck} implies
that $u_{\rm obs}(t,x)$ for small $t>0$ can only contain the following
Lax waves: a shock wave propagating from the second discontinuity
point $x=\xi_2$ towards $x=b$, a rarefaction wave propagating from the
first discontinuity point $x=\xi_1$ with positive speeds and,
possibly, a shock wave propagating from $x=\xi_1$ towards $x=a$. The
presence of the latter shock depends on the value $k_1$, which is
unknown to the observer. We can now proceed to the proof of the lemma.

\medskip

\n {\it ~(i)} Assume that at time $t=T_1$ a jump in $u_{\rm obs}(T_1,
\dott)$ appears at $x=a$. This jump corresponds to a shock $u$-wave
arriving from $x=\xi_1$ with negative speed. Let $\bar u_o=u_{\rm
  obs}(T_1,a-)$ and $v_o=u_{\rm obs}(T_1,a+)$ denote the densities
separated by the shock with propagation speed $\sigma_a$. Then we can
immediately deduce
\begin{equation*}
  k_1= k_o\,{f(v_o) \over f(u^m)}\,,
\end{equation*}
and thus
\begin{equation*}
  \xi_1 = a-\sigma_a T_1\,.
\end{equation*}
Since the shock cannot have interacted with any other wave, due to the
fact that both $\bar u_o$ and $k$ are constant in $[a,\xi_1)$, these
values represent the only possible choice of $k_1,\xi_1$ which
generates the observed shock.

Hence, to complete the reconstruction of $k(x)$ it only remains to
find $\xi_2$. This can be done by using the second condition
in~\eqref{eq:cond_nstat_1_thm}. At $t=T_2$, a new wave appears in
$x=b$ and it is a shock. Let $\sigma_b$ be the positive speed of this
shock and $v_1=u_{\rm obs}(T_1,b-)$ and $\bar u_o=u_{\rm obs}(T_1,b+)$
be the states separated by the shock. If
\begin{equation}\label{eq:verify}
  k_o\,{f(v_1) \over f(\bar u_o)}=k_1\,,
\end{equation}
then we can conclude that the shock has reached $x=b$ without
interacting with any other wave, and find $\xi_2$ as
\begin{equation*}
  \xi_2 = b-\sigma_b T_2\,.
\end{equation*}
Otherwise, if~\eqref{eq:verify} does not hold, $v_1$ is not the
original left state of the shock. Hence, the observed shock is in fact
the result of the interaction between the shock generated at $t=0+$ at
$x=\xi_2$ and the faster rarefaction wave generated at $t=0+$ at
$x=\xi_1$. In this case, to find $\xi_2$ we exploit the following
relation
\begin{equation}\label{eq:xi2}
  \xi_2-\xi_1 = \left( T_2 -\, {b-\xi_2 \over k_o f'(v_1)}\right)
  k_1 f'(w_1)\,, 
\end{equation}
where $w_1$ is such that $w_1<u^m$ and $k_of(v_1)=k_1 f(w_1)$. Observe
that~\eqref{eq:xi2} states that the rarefaction wave observed in
$(T_2,b)$ has traveled with speed $k_1 f'(w_1)$ in $(\xi_1,\xi_2)$ and
with speed $k_of'(v_1)$ in $(\xi_2,b)$. This completes the proof
when~\eqref{eq:cond_nstat_1_thm} holds.

\medskip

\n{\it (ii)} Set now
\begin{equation*}
  \tau:=\inf\left\{s\in~(T_1,T_2)~;~u(s,b)>u(T_1,b-)\right\}\,.
\end{equation*}
If $\tau>T_1$, then the shock wave generated at $x=\xi_2$ has reached
the observable region $[b,\infty)$ without interacting with the
centered rarefaction wave generated at $x=\xi_1$. Hence, we can repeat
the procedure followed in {\it (i)} and let $\bar u_o=u_{\rm
  obs}(T_1,b+)$, $v_o=u_{\rm obs}(T_1,b-)$ be the densities separated
by the shock which reaches $x=b$ at $t=T_1$ and $\sigma_b$ be its
positive propagation speed. Then,
\begin{equation*}
  k_1= k_o\,{f(v_o) \over f(\bar u_o)}\,,
\end{equation*}
and
\begin{equation*}
  \xi_2 = b-\sigma_b T_1\,.
\end{equation*}
To find $\xi_1$ we now use the fact that for $\tau<T_2$ the
rarefaction appears at $x=b$. Indeed, we know that the rarefaction
wave taking the value $v_o$ has traveled with speed $k_of'(v_o)$ in
$(\xi_2,b)$ and with speed $k_1 f'(w_o)$ in $(\xi_1,\xi_2)$, where
$w_o$ is such that $w_o<u^m$ and $k_of(v_o)=k_1 f(w_o)$. Hence, we can
exploit the relation
\begin{equation*}
  \xi_2-\xi_1 = \left( \tau -\, {b-\xi_2 \over k_o f'(v_o)}\right) 
  k_1 f'(w_o)\,,
\end{equation*}
which gives $\xi_1$.

On the other hand, if $\tau=T_1$, then the wave observed in $(\tau,b)$
is a rarefaction wave followed by an adjacent shock wave, and this
means that the interaction between the faster rarefaction wave and the
slower shock wave in $(\xi_2,b]$ has already begun. In this case, we
can still proceed as above: let $\bar u_o=u_{\rm obs}(T_1,b+)$,
$v_1=u_{\rm obs}(T_1,b-)$ be the densities separated by the shock and
$\sigma_b$ be its propagation speed and set
\begin{equation*}
  k_1= k_o\,{f(v_1) \over f(\bar u_o)}\,,
\end{equation*}
\begin{equation*}
  \xi_2 = b-\sigma_b T_1\,,
\end{equation*}
and
\begin{equation*}
  \xi_2-\xi_1 = \left( T_1 -\, {b-\xi_2 \over k_o f'(v_1)}\right)
  k_1 f'(w_1)\,,
\end{equation*}
where $w_1$ is again such that $w_1<u^m$ and $k_of(v_1)= k_1
f(w_1)$. These values of $(k_1,\xi_1,\xi_2)$ provide a solution
$u_{(k_1,\xi_1,\xi_2)}$ which coincides with $u_{\rm obs}$ outside
$(a,b)$, but they are not in general the only ones with such a
property.~~$\diamondsuit$

\section{Explicit reconstruction examples}

\begin{example}\label{ex:explicit1} Assume we want to reconstruct the flux in a sedimentation 
model where the local concentration $u$ of solids evolves according 
to~\eqref{eq:introhom} with $f$ of class $\con^2([0,1])$ and concave--convex 
with a single maximum point  $u_{\rm max}$ and a single inflection point 
$u_{\rm infl}$ (cf. for instance~\cite{BD} for examples of applications of this 
kind of models). Assuming that $f(0)=0$, we now apply the reconstruction 
procedure depicted in the proof of Theorem~\ref{thm:1} so to reconstruct a 
piecewise affine approximation of $f$ which coincides with the real flux in 
the points $u\in\{0,1/4,1/2,3/4,1\}$.
We first test the following Riemann initial data
$$
u_{o,1}(x)=\left\{
    \begin{array}{ll}
      0, & \,\,x<0,\\
      1/4, & \,\,x>0,
    \end{array}
  \right.
  \qquad
u_{o,2}(x)=\left\{
    \begin{array}{ll}
      1/4, & \,\,x<0,\\
      1/2, & \,\,x>0,
    \end{array}
  \right.
$$
$$
u_{o,3}(x)=\left\{
    \begin{array}{ll}
      1/2, & \,\,x<0,\\
      3/4, & \,\,x>0,
    \end{array}
  \right.
  \qquad
u_{o,4}(x)=\left\{
    \begin{array}{ll}
      3/4, & \,\,x<0,\\
      1, & \,\,x>0,
    \end{array}
  \right.
$$
and we observe the corresponding solutions $u_i$ ($i=1,\ldots,4$) at time $T=1$. 
Assume that the functions $u_1(1,\dott)$, $u_2(1,\dott)$ and $u_3(1,\dott)$ consist
of a single shock wave located respectively at $x_1=3$, $x_2=1$ and $x_3=-4/5$
and that $u_4(1,\dott)$ consists of a shock, located at $x_4=-5/2$ and separating 
the states $u_l=3/4$ and $u_m=7/8$, followed by an adjacent rarefaction wave, 
continuously increasing from $u_m=7/8$ to $u_r=1$ as $x\in [-5/2,0]$ 
(see Figure~\ref{fig:explicit_sol}). Also assume that the solution $u_4(1,\dott)$ satisfies 
along the rarefaction $\int_{-5/2}^0 u_4(1,x)\,dx=9/4$.

  \begin{figure}\begin{center}
%%      \psfrag{x}{${\scriptscriptstyle x}$} 
%%      \psfrag{z}{${\scriptscriptstyle 0}$} 
%%      \psfrag{a}{${\scriptstyle 0}$}
%%      \psfrag{b}{${\scriptstyle {1\over 4}}$}
%%      \psfrag{a1}{${\scriptstyle {1\over 4}}$}
%%      \psfrag{b1}{${\scriptstyle {1\over 2}}$}
%%      \psfrag{a2}{${\scriptstyle {1\over 2}}$}
%%      \psfrag{b2}{${\scriptstyle {3\over 4}}$}
%%      \psfrag{a3}{${\scriptstyle {3\over 4}}$}
%%      \psfrag{b3}{${\scriptstyle 1}$}
%%      \psfrag{k}{${\scriptstyle {7\over 8}}$} 
%%      \psfrag{u1}{${\scriptscriptstyle \!\!\!\!\!\!\!\!u_1(1,\dott)}$}
%%      \psfrag{u2}{${\scriptscriptstyle \!\!\!\!\!\!\!\!u_2(1,\dott)}$}
%%      \psfrag{u3}{${\scriptscriptstyle \!\!\!\!\!\!\!\!u_3(1,\dott)}$}
%%      \psfrag{u4}{${\scriptscriptstyle \!\!\!\!\!\!\!\!u_4(1,\dott)}$}
      \includegraphics[width=\textwidth]{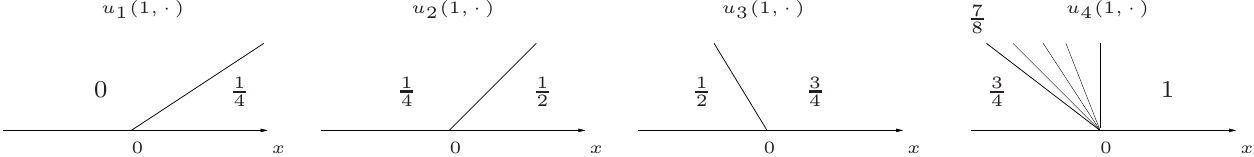}
      \caption{The solutions to the Riemann problems used in the reconstruction
      of Example~\ref{ex:explicit1}.\label{fig:explicit_sol}}
  \end{center}\end{figure}

\n This already allows to conclude that
$$
f\left({1\over 4}\right)=f(0)\,+\,3\,\cdot\,{1\over 4}\,=\,{3\over 4}\,,
\qquad\qquad
f\left({1\over 2}\right)=f\left({1\over 4}\right)\,+\,1\,\cdot\,{1\over 4}\,=\,1\,,
$$
$$
f\left({3\over 4}\right)=f\left({1\over 2}\right)\,-\,{4\over 5}\,\cdot\,{1\over 4}\,=\,{4\over 5}\,.
$$
Concerning the computation of $f(1)$, the presence of both a shock and a rarefaction wave
gives use additional knowledge. Indeed, from the presence of the shock wave separating 
$u_l$ and $u_m$ we deduce that
$$
f\left({7\over 8}\right)=f\left({3\over 4}\right)\,-\,{5\over 2}\,\cdot\,{1\over 8}\,=\,{7\over 16}\,,
$$
while from the presence of the rarefaction wave separating $u_m$ and $u_r$ we deduce
\begin{align*}
f(1)&=f\left({7\over 8}\right)\,+\int_{7/8}^1 f'(u)\,du\\
&=f\left({7\over 8}\right)\,+ \,1\,\cdot\,0\, - \,{7\over 8}\,\cdot\,\left(-\,{5\over 2}\right)\,
-\int_{-5/2}^0 u_4(1,x)\,dx=\,{3\over 8}\,,
\end{align*}
where we have applied Lemma~\ref{lem:inverse_int} in the second equality. Our approximate 
flux is then the affine function shown in Figure~\ref{fig:explicit_flux}--left, joining the points 
$(0,0)$, $(1/4,3/4)$, $(1/2,1)$, $(3/4,4/5)$, $(7/8,7/16)$, and $(1,3/8)$ and we know for sure 
that such an approximation coincides with the real flux in each of these points.

  \begin{figure}[b]\begin{center}
%%      \psfrag{x}{${\scriptstyle u}$} 
%%      \psfrag{z}{${\scriptstyle 0}$} 
%%      \psfrag{u}{${\scriptstyle 1}$} 
%%      \psfrag{1q}{${\scriptstyle {1\over 4}}$}
%%      \psfrag{2q}{${\scriptstyle {1\over 2}}$}
%%      \psfrag{3q}{${\scriptstyle {3\over 4}}$}
%%      \psfrag{fn}{${\scriptstyle ~~~~~~~~f_\nu(u)}$} 
%%      \psfrag{f}{${\scriptstyle ~~~~~~~~f(u)}$} 
%%
      \includegraphics[width=.99\textwidth]{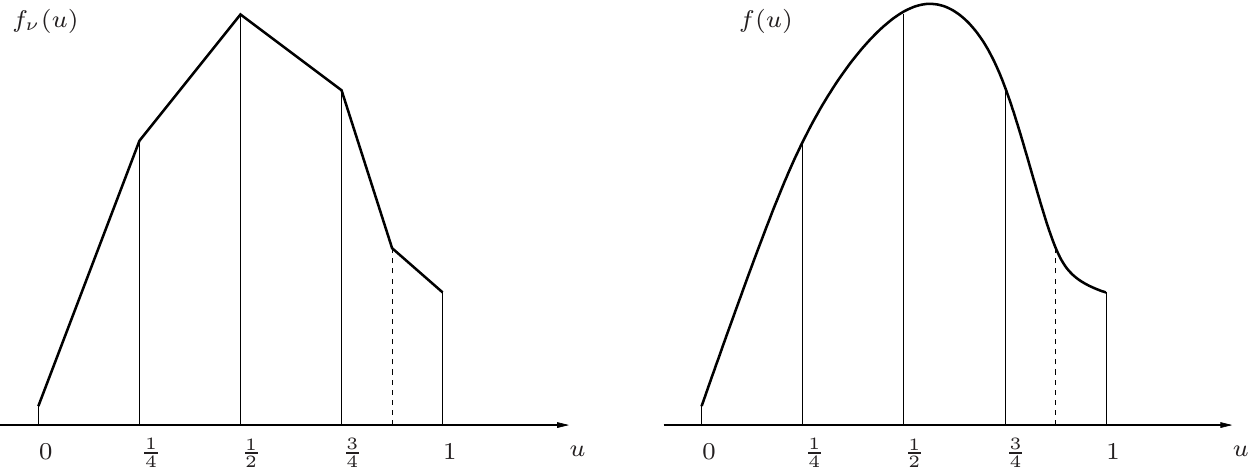}
      \caption{Left: Our reconstructed concave--convex flux $f_\nu$ in Example~\ref{ex:explicit1}. Right: 
      A possible real flux $f$ for the considered problem.\label{fig:explicit_flux}}
  \end{center}\end{figure}

\n It is now easy to deduce that  $u_{\rm max}$ must belong to either the interval $[1/4,1/2]$
or the interval $[1/2,3/4]$ and that  $u_{\rm infl}$ must belong to the interval $[3/4,7/8]$. If we want
a more precise localization of one of these point, it is enough to repeat the procedure splitting
the desired interval in suitably small subintervals and reconstructing the values of $f$ in 
intermediate points.
\end{example}

\begin{example}\label{ex:explicit2} We now consider an application of Theorem~\ref{thm:4} to the 
traffic flow model described in Example~\ref{ex:highway}. Namely, we consider a conservation
law~\eqref{eq:intrononh} describing the evolution of the car density $u\in [0,1]$ on a road, 
under the assumptions that $f(u)=u(1-u)$, that the initial distribution of cars is constant
$u(0,x)\equiv u_o$ and that we can only measure the solution $u(t,x)$ for $t\in [0,1]$ and 
$x\in\R\setminus (0,2)$. To fix the idea assume $u_o=1/3$.

\n Assume that at time $t=0$, the constant flux coefficient
$k_{\rm old}(x)\equiv 1$ is suddenly replaced by a piecewise
constant function $k(x)$ satisfying {\bf (H3)}, due to some car accidents 
occurring inside the unobservable interval $(0,2)$.
Since the initial data $u\equiv 1/3$ is not a stationary
solution for~\eqref{eq:intrononh} with discontinuous flux $k(x)f(u)$,
at time $t=0^+$ the solution will contain some traveling waves
around the discontinuity points for $k$.

\n Aim of our reconstruction procedure is to locate the position of the 
accidents by observing the density $u$ near
the ends of the unobservable region $x=0$ and $x=2$. In particular, if
\begin{itemize}
\item at time $T_1=0.5$ a shock wave emerges at $x=0$, separating the states
$u_o=1/3$ and $v_o=5/6$ and traveling with speed $\sigma_a=-1/6$,
\item at time $T_2=0.66$ a shock wave emerges at $x=2$, separating the states
$v_1=1/6$ and $u_o=1/3$ and traveling with speed $\sigma_b=1/2$, 
\end{itemize}
then we can apply the procedure depicted in the proof of Theorem~\ref{thm:4} 
(cf. Lemma~\ref{lem:4-2}), to conclude that the unique triple $(k_1,\xi_1,\xi_2)$ 
determining a flux $k$ of the form~\eqref{eq:disck} is given by
$$
k_1=\,{f(v_o)\over f(1/2)}\,=\,{5\over 9}\,,
\qquad
\xi_1=0 - \sigma_a \cdot T_1=\,{1\over 12}\,,
\qquad
\xi_2=2 - \sigma_b \cdot T_2=\,{167\over 100}\,.
$$
This means we can conclude that the accidents have created a partially obstructed region
$[{1\over 12}\,,\,{167\over 100}]\subset [0,2]$ where the flux is reduced to ${5\over 9}\, u(1-u)$.

\end{example}

\section{Conclusions}
In this paper we have presented some new results concerning inverse
problems for scalar conservation laws of the
form~\eqref{eq:intrononh}. Namely, for homogeneous
equations~\eqref{eq:introhom}, we have presented a reconstruction
procedure to find piecewise affine interpolations $f_\nu$ of any
piecewise $C^{1,1}$ flux $f$ having a finite number of inflection
points, under the unique assumption that solutions to Riemann problems
are observable at a fixed time $T>0$. No a priori assumption is
requested on the smoothness of the observed solution, or on its jumps
structure. The reconstructed flux is accurate in the following sense:
solutions to Cauchy problems for the conservation law with flux
$f_\nu$ are close in $\Lsp^1$ to the solutions for the conservation
law with exact flux $f$.

For general inhomogeneous equations~\eqref{eq:intrononh}, we have
first proved that being able to observe in $[0,T]\times\R$ the
solutions to Cauchy problems, for an arbitrarily small time $T$, is
sufficient to obtain a piecewise affine approximation of $f$, and the
precise form of $k(x)$ for $x$ in any compact interval $J\subseteq\R$.

Then, motivated by applications to traffic flow models, we have
studied the same inverse problem when the solutions are only
observable in part of the domain, due to the presence of some
inaccessible spatial region $I\subseteq\R$, and the goal is to
reconstruct $k(x)$ also inside $I$. In this case, even assuming the
observation of the solution for a long time interval $[0,T]$, the
function $k$ in the unobservable region can only be recovered under
the strict assumption that $k$ has no more than two jumps in
$I$. Unfortunately, this is not just a mathematical obstacle or a
limitation in the results we have presented: The examples in
Section~\ref{sec:counterex} show that if three or more jumps are
present, then in many situations we end up with an infinite number of
piecewise constant functions $k(x)$ which all give the same solution
in the observable region $\R\setminus I$.

In view of these examples, it seems clear that inverse problems for
inhomogeneous conservation laws, when only partial observability of
the solution is assumed, are in general ill-posed. Therefore, the next
steps in the study of inverse problems should focus the attention
either on specific inhomogeneous scalar models or on homogeneous
systems of hyperbolic conservation laws.  In the former case, one can
hope that physical features of the particular model considered help in
order to obtain well-posedness. In the latter case, one can try to
exploit the front-tracking algorithm to choose initial data which are
particularly well suited for the reconstruction, in the spirit of
Theorem~\ref{thm:1}. This problem is much more difficult for systems
than for the single equation, due to the possibly complicate wave
structure of the solutions, but some positive result could be
possible, at least in the case of Temple class systems, which have
coinciding shock and centered rarefaction waves. Alternatively, one
could try to adapt the least square method, used
in~\cite{JamesSepulvedaSIAM}, to the generalized differentiability
structures which have been introduced for systems of conservation laws
by Bressan and Marson~\cite{BressanMarson}, and look for a
reconstruction of the flux as the minimizer of a cost
like~\eqref{eq:JScost}. Also with this approach, however, the case of
systems is substantially more difficult than the case of the single
equation, and it is not clear which regularity can be expected from
the cost functional.

\medskip

\n{\bf Acknowledgements.} We thank the anonymous referees, who read 
very carefully this paper,  for their comments and remarks which helped to 
improve the exposition.

\appendix
\section*{Appendix}
\renewcommand{\thesection}{A}

Here we collect some auxiliary results which have been used in the
paper. We start from a simple property from standard calculus which
have been exploited in the proof of Theorem~\ref{thm:1}.
\begin{lemma}\label{lem:inverse_int} If $\gamma\colon [a,b]\to\R$ is a
  continuous and strictly monotone function, then
  \begin{equation}\label{eq:areaformula}
    \int_{\gamma(a)}^{\gamma(b)} \gamma^{-1}(s) \, ds = \gamma(b)b - \gamma(a)a
    - \int_a^b \gamma(t)\,dt\,.
  \end{equation}
\end{lemma}
% \n{\it Proof.} Assume that $\gamma$ is strictly increasing, being
% the strict decreasing case analogous. Then, by Fubini's Theorem, we
% have
% $$
% \int_a^b \gamma(t)\,dt= \int_a^b \int_0^{\gamma(t)}\,ds\,dt=
% \gamma(a)(b-a)+ \int_{\gamma(a)}^{\gamma(b)}
% \int_{\gamma^{-1}(s)}^b\,dt\,ds
% $$
% Hence,
% \begin{align*}
%   \int_a^b \gamma(t)\,dt &+ \int_{\gamma(a)}^{\gamma(b)}
%   \gamma^{-1}(s)\,ds\\ &~= \gamma(a)(b-a)+
%   \int_{\gamma(a)}^{\gamma(b)}\left( \int_{\gamma^{-1}(s)}^b\,dt +\gamma^{-1}(s) \right)\,ds\\
%   &~=\gamma(a)(b-a)+ b \left(\gamma(b)-\gamma(a)\right)\\
%   &~=b\,\gamma(b)- a\, \gamma(a)\,,
% \end{align*}
% which is exactly~\eqref{eq:areaformula}.~~$\diamondsuit$

% \medskip

Next, we recall a well-posedness result for scalar conservation laws
that we have exploited to show that~\eqref{eq:recest}
implies~\eqref{eq:close}. The proof can be found in~\cite[Theorem
2.3]{HoldenRisebroBook}.

\begin{theorem}\label{thm:lip_est} Let $f,g$ be Lipschitz continuous
  functions, and assume $\hat u\in\BV(\R)$.  Denote by $u^f$ and
  $u^g$, respectively, the solutions to the Cauchy problems for
$$
\partial_t u + \partial_x f(u) =0\,, \qquad\qquad
\partial_t u + \partial_x g(u) =0\,,
$$
with initial data $u(0,x)=\hat u(x)$. Then, there exists a positive
constant $C$ such that for all $T\geq 0$
$$
\| u^f(T,\dott) - u^g(T,\dott) \|_{\Lsp^1(\R)} \le CT
\mathrm{\Lip}(f-g)\,.
$$
\end{theorem}
Finally, we offer a brief description of the Riemann solver
for~\eqref{eq:intrononh} defined and studied
in~\cite{KlingRisebro,Risebrointro}, which we have used extensively in
the proofs of Theorems~\ref{thm:2}--\ref{thm:4}. First of all, taking
flux functions $k(x)$ and $f(u)$ that satisfy {\bf (H1)}, we observe
that the inhomogeneous equation~\eqref{eq:intrononh} can be studied by
considering an auxiliary system of conservation laws
\begin{equation}\label{eq:aux_sys}
  \left\{\begin{array}{l}
      \partial_t u + \partial_x (kf(u)) =0,\\
      \partial_t k = 0,
    \end{array}
  \right.
\end{equation}
which represents the conservation of the quantity $v=(u,k)$ with flux
$g(v)=(kf(u),0)$. The aim of this auxiliary system $\partial_t v
+ \partial_x g(v) =0$ is to help in the study of the behavior of the
solution to~\eqref{eq:intrononh} at discontinuities of
$k(x)$. However, such an auxiliary system is non-strictly hyperbolic,
since waves of the second family (i.e., related to the second
equation) all have null speed, while waves of the first family (i.e.,
related to the original scalar equation) can have positive or negative
speeds depending on the sign of $f'$. Hence, the
system~\eqref{eq:aux_sys} requires some additional attention.

% Systems of this kind are sometimes called ``resonant systems'' and
% have been the subject of a wide literature in the past decades.
%
The properties we need to know here are the following:
\begin{itemize}\item
  a solution to any Riemann problem for~\eqref{eq:aux_sys} can be
  constructed by following~\cite{GimseRisebroProc};

\item the solution is unique, provided an additional ``entropy''
  condition holds at jumps of $k$ (the precise condition will be
  discussed below);

\item the construction allows one to build a converging front-tracking
  approximation for general Cauchy problems
  following~\cite{KlingRisebro}.
\end{itemize}
Namely, the construction proceeds as follows. The Rankine--Hugoniot
conditions for~\eqref{eq:aux_sys} can be written as
$$
kf(u)-k'f(u')=\lambda (u-u')\,, \qquad \qquad 0=\lambda (k-k')
$$
for a discontinuity separating states $(u,k)$,$(u',k')$ and traveling
with speed $\lambda$. In other words, either $k=k'$ and we have a
discontinuity in $u$ only, or $\lambda=0$ and the states separated by
the stationary jump satisfy
\begin{equation}\label{eq:rh_stationary}
  kf(u)=k'f(u')\,.
\end{equation}
In particular, all discontinuities in $k$ give origin to a stationary
jump in the solution, with both $k$ and $u$ being discontinuous across
the jump.

\n Under the assumption {\bf (H1)}, for a fixed state $u$ and fixed
constants $k,k'$ in general there exist two solutions $u'=v_1$ and
$u'=v_2$ to~\eqref{eq:rh_stationary}, and they satisfy $v_1\le
u^m\le v_2$. The admissibility condition (or ``entropy'' condition)
mentioned above, which is needed to select a single state $u'$ at the
stationary jumps, is the following: the admissible state is the one
which realizes $\min\{|u-v_1|,|u-v_2|\}$. In general, fix a jump of
$k(x)$ and denote the states adjacent to the discontinuity and
satisfying~\eqref{eq:rh_stationary} by $k^-,u^-$ and $k^+,u^+$. Then
the solution connecting these states is \emph{entropy admissible} if
and only if $u^+,u^-$ satisfy
\begin{equation}\label{eq:smallest_jump}
  |u^+-u^-|=\min\big\{|v-v'|~;~k^+f(v)=k^-f(v')\big\}\,,
\end{equation}
i.e., if they minimize the quantity $|v-v'|$ among all pairs
satisfying~\eqref{eq:rh_stationary}. In~\cite{GimseRisebroProc} it was
shown that this ``smallest jump'' condition is equivalent to a viscous
profile entropy condition for the auxiliary system~\eqref{eq:aux_sys},
justifying the use of the word entropy also in the context
of~\eqref{eq:intrononh}.

\medskip

%% As in most of the result concerning this kind of systems of
%% conservation laws, we will solve the problem in a different system
%% of coordinates. Let the singular function $\Psi$ be defined as
%% follows
%%%
%% \begin{equation}\label{eq:psi}
%%\Psi(u,k):= 
%%\left\{\begin{array}{ll}
%%k\,{u-u^m \over |u-u^m|} \left( 1-{f(u) \over f(u^m)} \right)  & \mbox{if } u\neq u^m\\
%%0 & \mbox{if } u=u^m
%%\end{array}\right.
%%\end{equation}
%%%
%%Then, we will consider the problem in the new coordinates $(z,k):=
%% \left(\Psi(u,k), k\right)$. The following Proposition collects
%% various properties of $\Psi$ that we will need.
%%
%% \vskip .2cm
%%
%% \n {\bf Proposition.}  Let $\Psi$ be the function defined
%% in~\eqref{eq:psi}. Then the following holds:
%%
%% \begin{description}
%% \item{(i)} $\Psi$ is continuous in its whole domain and smooth for
%%   $u\neq u^m$.
%% \item{(ii)} For any fixed $u\in [0,1]$, the function
%%   $\Psi(u,\dott)$ satisfies $\forall k_1,k_2\in~]0,\bar k]$
%%$$
%%|\Psi(u,k_1)-\Psi(u, k_2)|\le |k_1-k_2|
%%$$
%% \item{(iii)} For any fixed $k\in ]0,\bar k]$, the function
%%   $\Psi(\dott,k)$ is strictly monotone increasing, maps $[0,1]$
%%   into $[-k,k]$ and satisfies $\forall u_1,u_2\in[0,1]$
%%$$
%%|\Psi(u_1,k)-\Psi(u_2, k)|\le k\,{||f'||_\infty \over
%%  ||f||_\infty}|u_1-u_2|
%%$$
%% \end{description}
%%
%%
%% \medskip
%%
\n Now assume that $k$ has a discontinuity at $x=0$ and that we are
given a Riemann initial data
\begin{equation}\label{eq:aux_data}
  v_o(x)=\begin{cases}
    (u^\ell,k^\ell), & \mbox{ if }x<0,\\
    (u^r,k^r), & \mbox{ if }x>0,
  \end{cases}
\end{equation}
for suitable constants $u^\ell,u^r\in[u_1,u_2]$ and
$k^\ell,k^r>0$. Following the conventions
of~\cite{GimseRisebroProc,KlingRisebro}, we will call $u$-waves the
Lax waves with constant $k$ (equivalently, the Lax waves of the first
family for~\eqref{eq:aux_sys}) and $k$-waves the Lax stationary waves
where $k$ changes and~\eqref{eq:rh_stationary} holds (equivalently,
the Lax waves of the second family for~\eqref{eq:aux_sys}). Then, the
solution to~\eqref{eq:aux_sys}--\eqref{eq:aux_data} can be constructed
as follows:

\medskip

\n{\bf Case 1.} Assume $k^\ell<k^r$ and $u^\ell\le u^m$. Then, if
$u^r\le u^m$ or $k^\ell f(u^\ell)<k^r f(u^r)$, the solution is given
by a stationary $k$-wave between $(k^\ell,u^\ell)$ and $(k^r,v)$, with
$v$ satisfying $v<u^m$ and $k^\ell f(u^\ell)=k^rf(v)$, followed by a
$u$-shock or a centered $u$-rarefaction with positive speed, between
the states $(k^r,v),(k^r,u^r)$. On the other hand, if $u^r>u^m$ and
$k^\ell f(u^\ell)\geq k^r f(u^r)$, the solution is given by a
$u$-shock, traveling with negative speed, between the states
$(k^\ell,u^\ell)$ and $(k^\ell,w)$, with $w>u^m$ and $k^\ell
f(w)=k^rf(u^r)$, followed by a $k$-wave separating $(k^\ell,w)$ and
$(k^r,u^r)$.
% the case $k^\ell f(u^\ell)= k^r f(u^r)$ is a limit case, in the
% sense that the shock between $w=u^\ell$ and $v$ has zero speed
% (because $f(w)=f(u^\ell)$) and ``disappears'' into the jump in $k$
% (it's easier to see the solution in the alternate coordinates used
% by~\cite{KlingRisebro}).

\medskip

\n{\bf Case 2.} Assume $k^\ell<k^r$ and $u^\ell> u^m$. Then, if
$u^r\le u^m$ or $k^r f(u^r)>k^\ell f(u^m)$, the solution is given by
a centered $u$-rarefaction with positive speed, between the states
$(k^\ell,u^\ell),(k^\ell,u^m)$, followed by a $k$-wave between
$(k^\ell,u^m)$ and $(k^r,v')$, with $v'$ satisfying $v'<u^m$ and
$k^\ell f(u^m)=k^rf(v')$, followed by a $u$-shock or a $u$-rarefaction
with positive speed, separating the states $(k^r,v'),(k^r,u^r)$. On
the other hand, if $u^r>u^m$ and $k^r f(u^r)\le k^\ell f(u^m)$, the
solution is given by a $u$-shock or a $u$-rarefaction with positive
speed, between the states $(k^\ell,u^\ell)$ and $(k^\ell,w')$, with
$w'\geq u^m$ and $k^\ell f(w')=k^rf(u^r)$, followed by a $k$-wave
separating $(k^\ell,w')$ and $(k^r,u^r)$.

\medskip

\n{\bf Case 3.} Assume $k^\ell>k^r$ and $u^r\le u^m$. Then, if
$u^\ell\geq u^m$ or $k^\ell f(u^\ell)> k^r f(u^m)$, the solution is
given by a $u$-shock or a centered $u$-rarefaction with negative
speed, separating the states $(k^\ell,u^\ell)$ and $(k^ \ell,v'')$,
with $v''$ satisfying $v''>u^m$ and $k^r f(u^m)=k^\ell f(v'')$,
followed by a $k$-wave between $(k^\ell,v'')$ and $(k^r,u^m)$,
followed by a $u$-rarefaction with positive speed, between the states
$(k^r,u^m),(k^r,u^r)$. On the other hand, if $u^\ell<u^m$ and $k^\ell
f(u^\ell)\le k^r f(u^m)$, the solution is given by a $k$-wave between
the states $(k^\ell,u^\ell)$ and $(k^r,w'')$, with $w''<u^m$ and $k^r
f(w'')=k^\ell f(u^\ell)$, followed by a $u$-shock or a $u$-rarefaction
with positive speed, separating the states $(k^r,w''),(k^r,u^r)$.

\medskip

\n{\bf Case 4.} Assume $k^\ell>k^r$ and $u^r> u^m$. Then, if
$u^\ell\geq u^m$ or $k^\ell f(u^\ell)\geq k^r f(u^r)$, the solution is
given by a $u$-shock or a centered $u$-rarefaction with negative
speed, between the states $(k^\ell,u^\ell)$ and $(k^ \ell,v''')$, with
$v'''$ satisfying $v'''>u^m$ and $k^r f(u^r)=k^\ell f(v''')$, followed
by a stationary $u$-wave separating $(k^\ell,v''')$ and $(k^r,u^r)$.
% the case $k^\ell f(u^\ell)= k^r f(u^r)$ is a limit case, in the
% sense that the shock between $v'''$ and $u^\ell$ has zero speed
% (because $f(v''')=f(u^\ell)$) and ``disappears'' into the jump in
% $k$ (it's easier to see the solution in the alternate coordinates
% used by~\cite{KlingRisebro}).
On the other hand, if $u^\ell<u^m$ and $k^\ell f(u^\ell)< k^r f(u^r)$,
the solution is given by a $k$-wave between the states
$(k^\ell,u^\ell)$ and $(k^r,w''')$, with $w'''<u^m$ and $k^r
f(w''')=k^\ell f(u^\ell)$, followed by a $u$-shock with positive
speed, between the states $(k^r,w'''),(k^r,u^r)$.

\medskip

This construction of a solution to the Riemann
problem~\eqref{eq:aux_sys}--\eqref{eq:aux_data} provides a Riemann
solver which allows to solve the Cauchy problem too, by means of a
standard front-tracking algorithm~\cite{BressanBook,
  HoldenRisebroBook}. The precise proof of the next theorem, and in
particular of the compactness of the approximation which allows to
apply Helly's theorem, can be found in~\cite{KlingRisebro}.

\medskip

\begin{theorem}\label{thm:KR} Let $f,k$ be flux functions satisfying
  {\bf (H1)}. Then, for every initial data $\hat u\in\BV(\R)$, the
  Cauchy problem for
$$
\partial_t u + \partial_x \big(k(x)f(u)\big)=0\,,
$$
with initial data $u(0,x)=\hat u(x)$ admits a weak solution $u$ such
that, for every time $t\ge 0$, $u(t,\dott)\in \Lsp^1(\R)$ is obtained
as the uniform $\Lsp^1_{\rm loc}$ limit of a front-tracking
approximation $u^\delta(t,\dott)$, constructed using the Riemann
solver described above.
\end{theorem}

\medskip

%% original Theorem required $\hat u$ in BV with respect to the other
%% variables $(\Psi,k)$, but in our settings this follows from $k$ and
%% $\hat u$ in BV of $\R$. indeed
%%$$
%%\sum_j |\Psi(x_i)-\Psi(x_j)|+|k(x_i)-k(x_j)|\le ||k||_\infty
%% ||f'||_\infty \TV(u_o) +2\TV(k)
%%$$
%%

\end{document}